\documentclass[12pt,twoside]{article}
\makeatletter
\let\@fnsymbol\@arabic
\makeatother

\usepackage{graphicx}
\usepackage{graphics}
\usepackage{amsmath,amssymb,amsthm}
\usepackage{authblk}

\theoremstyle{plain}
\newtheorem{thm}{Theorem}[section]

\newtheorem{lem}{Lemma}[section]

\theoremstyle{remark}

\newtheorem*{pf}{Proof}

\def\E{{\mathbb {E}}}                       \def\D{{\mathbb {D}}}

\def\FD{{\mathcal F}}

\def\o{{\rm {o}}} 
\def\O{{\rm {O}}}

\date{\vspace{-5ex}}
\title{The $N$-stars network evolution model}
\author[*]{{\sc Istv\'an Fazekas}\thanks{email: fazekas.istvan@inf.unideb.hu}} 
\author[*]{{\sc Csaba Nosz\'aly}\thanks{email: noszaly.csaba@inf.unideb.hu}} 
\author[*]{{\sc Attila Perecs\'enyi}\thanks{email: perecsenyi.attila@inf.unideb.hu}}
\affil[*]{Faculty of Informatics, University of Debrecen, P.O. Box 12\\ 4010 Debrecen, Hungary}

\begin{document}
\maketitle
\begin{abstract}
A new network evolution model is introduced in this paper.
The model is based on co-operations of $N$ units.
The units are the nodes of the network and the co-operations are indicated by directed links.
At each evolution step $N$ units co-operate which formally means that they form a directed $N$-star subgraph.
At each step either a new unit joins to the network and it co-operates with $N-1$ old units or
$N$ old units co-operate.
During the evolution both preferential attachment and uniform choice are applied.
Asymptotic power law distributions are obtained both for the in-degrees and the out-degrees.
\end{abstract}
\renewcommand{\thefootnote}{}
{\footnotetext{


\textsc{Attila Perecs\'enyi was supported through the New National Excellence Program of the Hungarian Ministry of Human Capacities.
Istv\'an Fazekas was supported by the construction EFOP-3.6.3-VEKOP-16-2017-00002;
the project was supported by the European Union, co-financed by the European Social Fund.}
}}

\section{Introduction}



Network science emerged during the past two decades (see \cite{barabasiBook}, \cite{durrett}).
It studies general features of real-world networks.
Such networks are the WWW, the Internet, the power grid, biological, social and trade networks.
In the Introduction of \cite{barabasiBook} A. L. Barab\'asi writes that
'A key discovery of network science is that the architecture of networks emerging in various domains of science, nature, and technology are similar to each other...'
and
'...we will never understand complex systems unless we develop a deep understanding of the networks behind them.'

In their seminal paper \cite{barabasi} Barab\'asi and Albert list various scale free large networks (actor collaboration, WWW, power grid, etc.), 
and they describe the preferential attachment model moreover, give an argument and simulation evidence that the preferential attachment rule leads to a scale-free network.
A network is called scale-free if its degree distribution is asymptotically power law, that is
$p_k \sim Ck^{-\gamma}$ as $k\to\infty$, where $p_k$ is the probability that a node has degree $k$.
Here and in what follows $a_k \sim b_k$ means that $\lim_{k\to\infty} a_k/b_k =1$.
The preferential attachment network evolution model is the following.
At every time step $t=2,3,\dots$ a new vertex with $N$ edges is added to the existing graph so that the edges link the new vertex to $N$ old vertices.
The probability $\pi_i$ that the new vertex will be connected to the old vertex $i$ depends on the degree $d_i$ of vertex $i$, so that $\pi_i= d_i/\sum_j d_j$,
where $\sum_j d_j$ is the cumulated sum of degrees.

We have to mention that long before the publication of \cite{barabasi}, Yule proposed a model for evolution of species where preferential attachment was present 
(see \cite{Yule}, see also \cite{Simon}).
Moreover, the precise mathematical formulation of the preferential attachment network evolution model  and a rigorous proof of 
the power law degree distribution in the preferential attachment model was given in Bollob\'as et al. \cite{bollobas} (see also \cite{durrett}, \cite{hofstad} and \cite{Chung-Lu}).
Nevertheless, in \cite{barabasi} Barab\'asi and Albert revealed the connection of preferential attachment and power law.
In his monograph \cite{hofstad} van der Hofstad underlines this connection 'A possible and convincing explanation for the occurrence of power-law degree sequences
is offered by the preferential attachment paradigm.'

The concept of preferential attachment and the scale-free property incited enormous research activity.
In connection with the mathematical models we also have to mention that the classical Erd\H{o}s-R\'enyi graph (see \cite{erRe1}, \cite{erRe2}, see also \cite{gil}) is not scale free.
Therefore new mathematical models were necessary to describe real-life networks.
For the mathematical theory see the monograph \cite{hofstad} written by van der Hofstad (see also \cite{durrett} and \cite{Chung-Lu}).
Concerning the general aspects of network theory one can consult the comprehensive book \cite{barabasiBook} by A. L. Barab\'asi.
In \cite{barabasiBook} a complete chapter is devoted to the scale free property. 
Based on previous studies of large real-life networks (WWW, Internet, e-mail, citation,...)
the author claims that scale-free property is a 'universal network characteristic'.
On the other hand he mentions that 'The ubiquity of the scale-free property does not mean that all real networks are scale-free.', and he lists networks not sharing this property.
In the literature, there are lot of papers devoted to the study of scale free property, but there are also papers not supporting this property.
For example the authors of the paper \cite{broido} claim that 'scale-free networks are rare' and 
'real-world networks exhibit a rich structural diversity that will likely require new ideas and mechanisms to explain'.
In  our paper we do not study any specific real-life network, but we offer a new mathematical model to build a network.
Our proposal is based on the star-like substructures of networks which on the one hand lead to a mathematically tractable model 
and on the other hand they seem to be plausible ingredients of real-life networks.

There are several versions of the preferential attachment model, here we can mention only a few of them.
In \cite{cooper} Cooper and Frieze introduced the following general graph evolution model.
At each step either a new vertex or an old one generates new edges.
In both cases the terminal vertices can be chosen either uniformly or according to the degrees of the vertices.
In \cite{osroumova} a general preferential attachment model (so called PA-class) is defined. 
Several known models (the LCD-model of \cite{bollobasRiordan}, the Holme-Kim model \cite{holmeKim}, the random Apollonian network \cite{Zhou}, 
the Buckley-Osthus-M\'ori model \cite{buckleyOsthus}, \cite{Mori}) belong to the PA-class.
In  \cite{osroumova} power law degree distribution was proved for the PA-class.
In \cite{Faz} the PA-class was extended to describe the evolution of certain populations.
In \cite{BaMo1}, \cite{BaMo2} and \cite{FIPB2} the above mentioned ideas of Cooper and Frieze \cite{cooper} were applied, 
but instead of the original preferential attachment rule, the vertices were chosen according to the weights of certain cliques.
An $N$-clique is a complete graph on $N$ vertices.
A clique can be considered as a particular model of a team.
That is any two members of the team are connected to each other.
However, there are other structures of cooperation among team members.

In this paper we shall consider star-like structures.
It means that there is a head of the team and all other members are connected to him/her.
E.g. a given person and his/her friends form a team, that is the center is the given person and the peripheral members are his/her friends.
Usually a person can play both roles.
E.g. John is the center in the team of his friends, but he is a peripheral member in the team of Peter's friends (assuming that John and Peter are friends).
This kind of double roles will be allowed in our model.
Examples of star-like structures can be found at companies, authorities, universities, etc.
Star topology is a usual structure in computer networks, see e.g. \cite{comp}.
Our model was motivated by star-like structures in the society and in technology.
However, our aim was to give a real mathematical model and not to describe a particular network.

Here we explain the evolution of our network in terms of persons.
The basic unit of the cooperation is a team of $N$ persons so that one of them plays central role and the others join to him/her.
So the structure of a team looks like a star on $N$ vertices.
A star on $N$ vertices (in short $N$-star) consists of a central vertex and $N-1$ peripheral vertices which are  connected to the central vertex.
Here we consider a star as a directed graph, the starting point of an edge is always a peripheral vertex and the target is the central vertex.
In our model the cooperation of $N$ persons always means that they form an $N$-star.
The cooperation of the same persons can be activated several times.
We allow multiple edges in order to show repeated cooperation. 
So we indicate the new cooperation by creating new edges.
For example, if our network consists of one $N$-star which was activated two times, 
then it has $2(N-1)$ directed edges so that its central vertex has in-degree $2(N-1)$ and each of the $N-1$ peripheral vertices has out-degree $2$. 

In our model the teams compete each other.
The strength of a team is measured by its weight.
If a team is activated again, then its weight is increased by 1.
The higher the weight of a team, the higher the chance that it will be activated again.
During the evolution new members can join to the network.
A newcomer has two possibilities.
Either he/she joins to an existing team or he/she creates a new team.
In the first case the newcomer chooses one of the existing teams according to the weights of the teams.
In the second case the newcomer chooses $N-1$ persons uniformly at random and the newcomer himself/herself will be the head of the new team. 
There are also evolution steps when there is no new person to join to the network.
In this case either an existing team is activated again or $N$ randomly chosen members of the network form a new team. 
The in-degree $d_1$ and the out-degree $d_2$ will describe the role of the person. 
$d_1/(N-1)$ is the number of cases when he/she was head of any team while
$d_2$ is the number of cases when he/she was a non-head member of any team.

In Section \ref{Model}, the precise mathematical description of the model is given.
Then scale-free property is proved both for in-degrees and out-degrees, see Theorem \ref{main}.
The proofs are presented in Section \ref{proofs}.
In the proofs the main probabilistic tools are the Doob-Meyer decomposition and convergence theorems for submartingales.
These are usual tools to obtain asymptotic results for random graphs (see, e.g., \cite{BaMo1}, \cite{BaMo2}, \cite{FIPB2}).
However, the models in \cite{BaMo1}, \cite{BaMo2}, \cite{FIPB2} were 'homogeneous' while in this paper we distinguish central and peripheral vertices.
Therefore our formulae are more complicated than those of the above mentioned papers.
So we could not use directly any calculation of the previous papers.
In this paper the challenge of the proofs was to handle multiple sequences and to guess the formulae to be proved by induction.

\section{The model and the main result} \label{Model}
\begin{figure}[!ht]
\centering
\includegraphics[width=5cm]{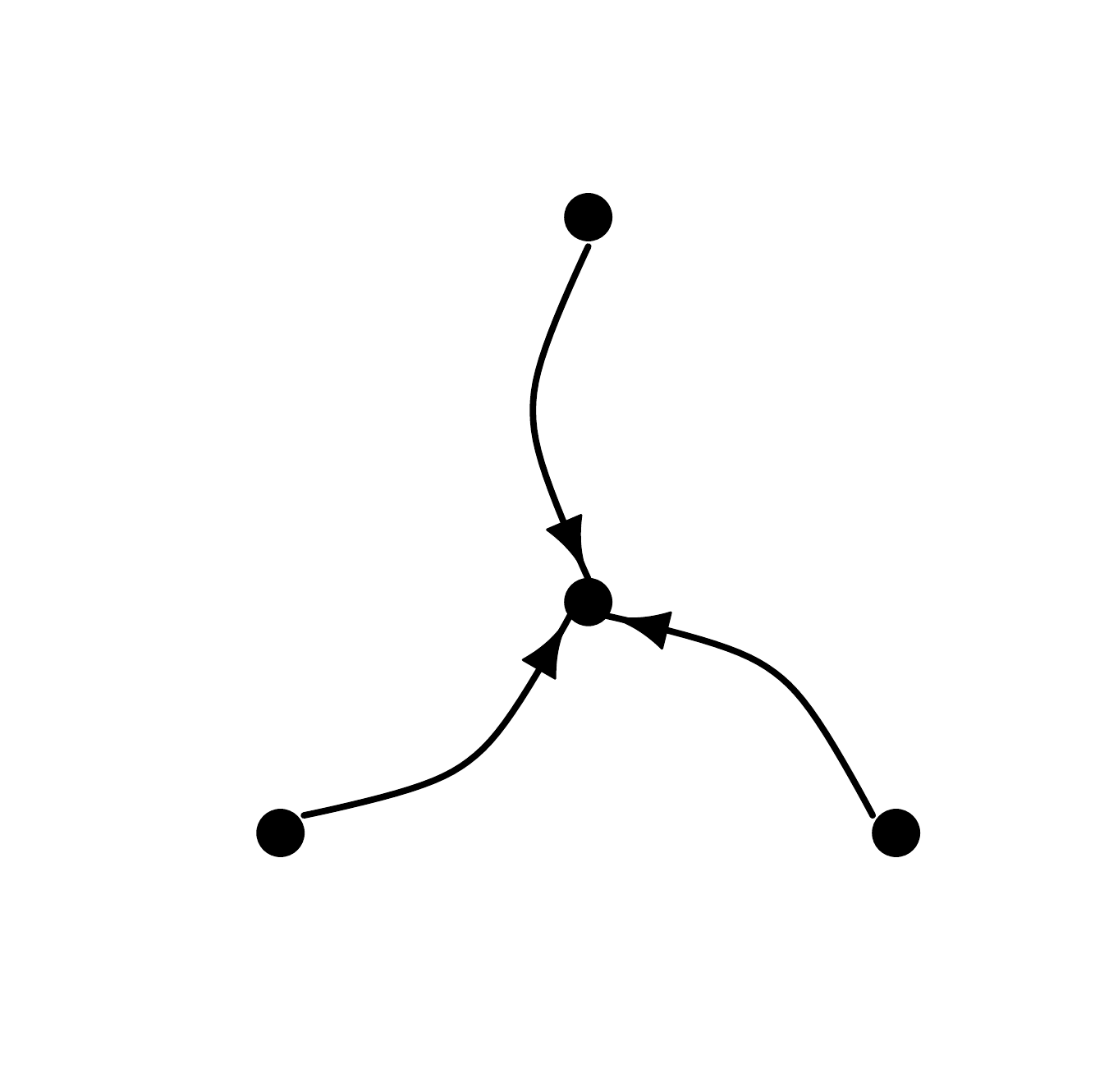}
\caption{The initial graph, $N=4$}\label{init}
\end{figure}

\begin{figure}[!ht]
  \centering
  \begin{minipage}[b]{0.4\textwidth}
    \includegraphics[width=5cm]{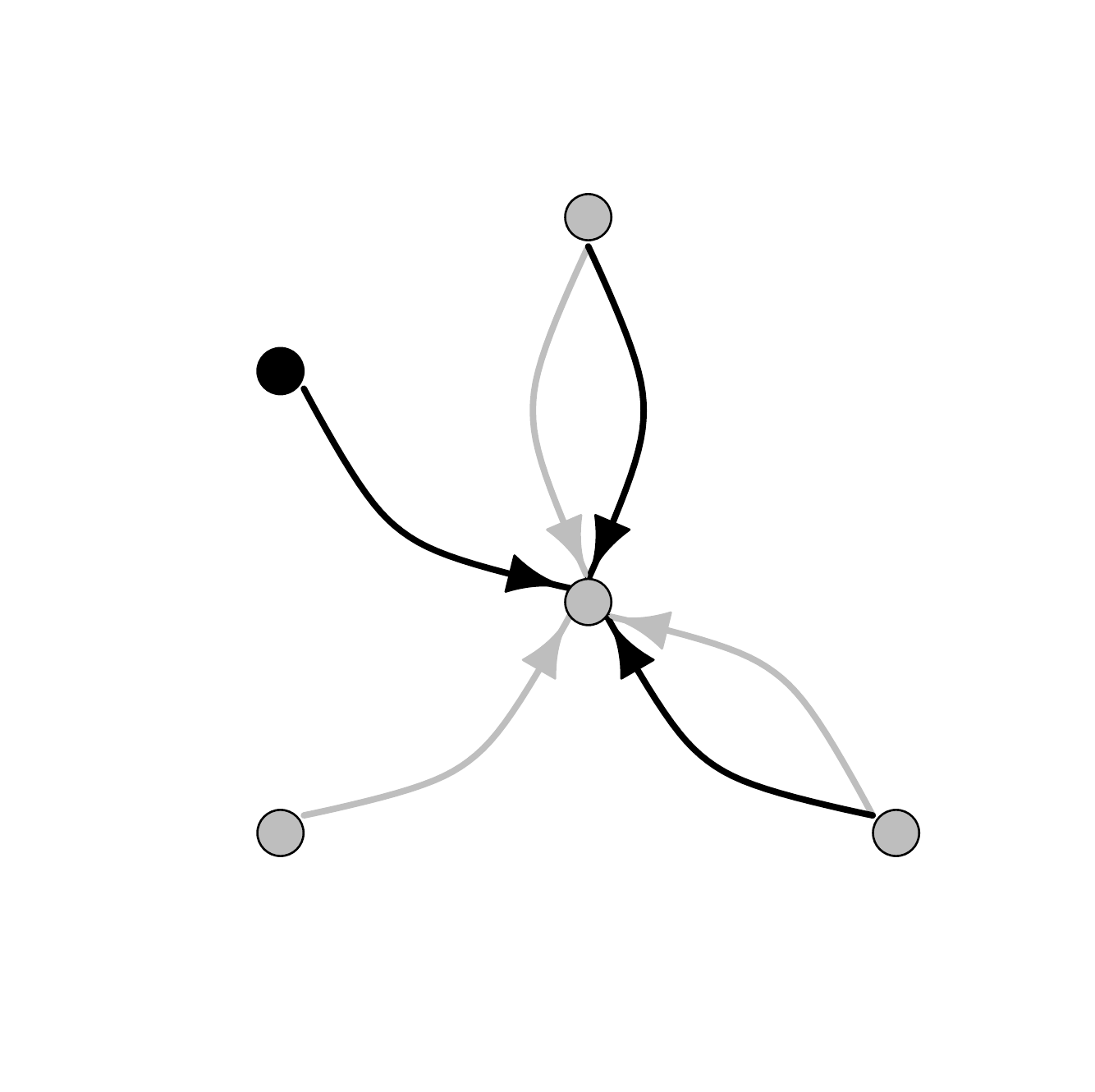}
    \caption{Case I/1 in the first step, $N=4$}\label{pr}
  \end{minipage}
  \hfill
  \begin{minipage}[b]{0.4\textwidth}
    \includegraphics[width=5cm]{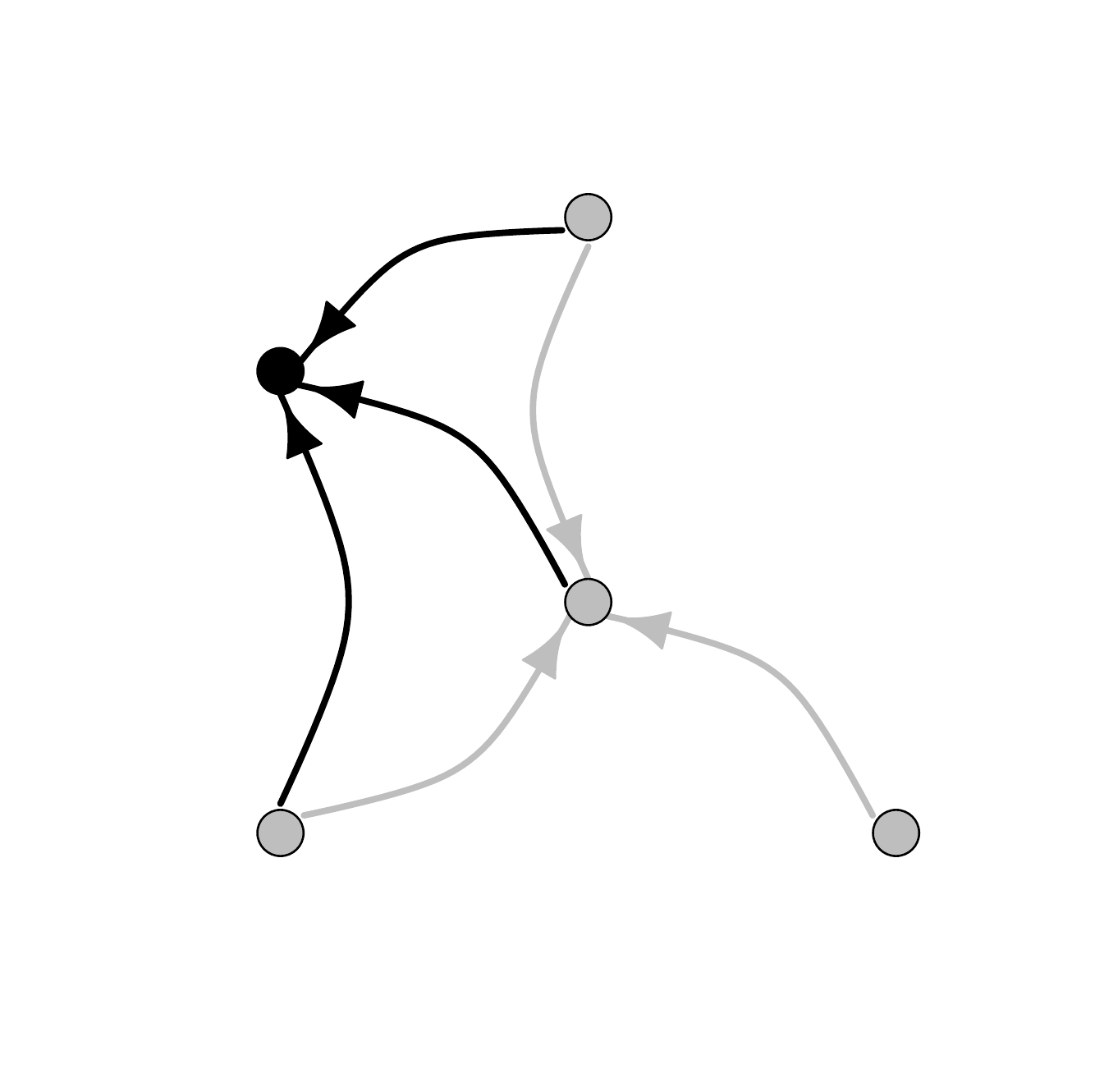}
    \caption{Case I/2 in the first step, $N=4$}\label{p1r}
  \end{minipage}
\end{figure}

\begin{figure}[!ht]
  \centering
  \begin{minipage}[b]{0.4\textwidth}
   \includegraphics[width=5cm]{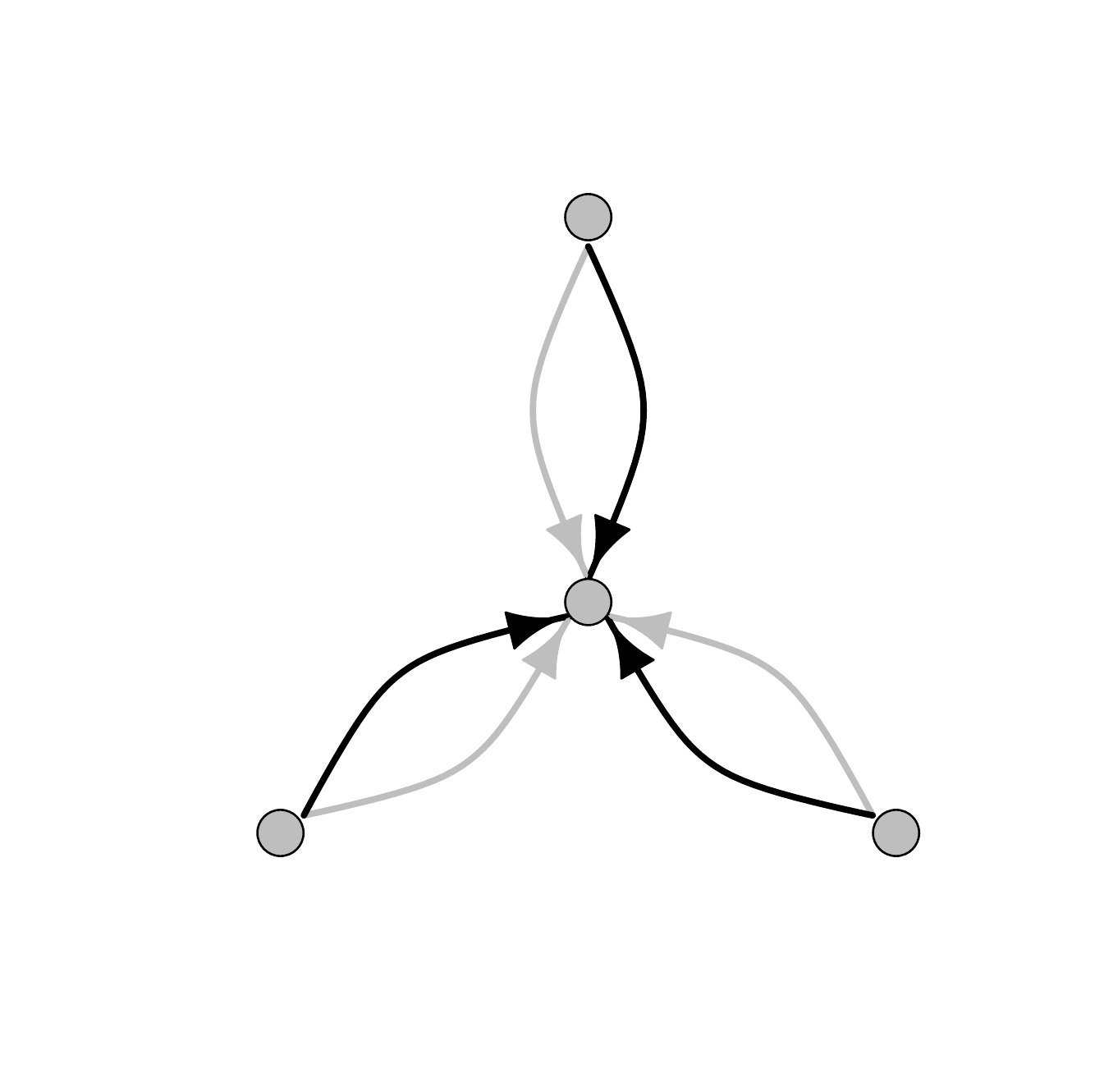}
    \caption{Case II/1 in the first step, $N=4$}\label{1pq}
  \end{minipage}
  \hfill
  \begin{minipage}[b]{0.4\textwidth}
   \includegraphics[width=5cm]{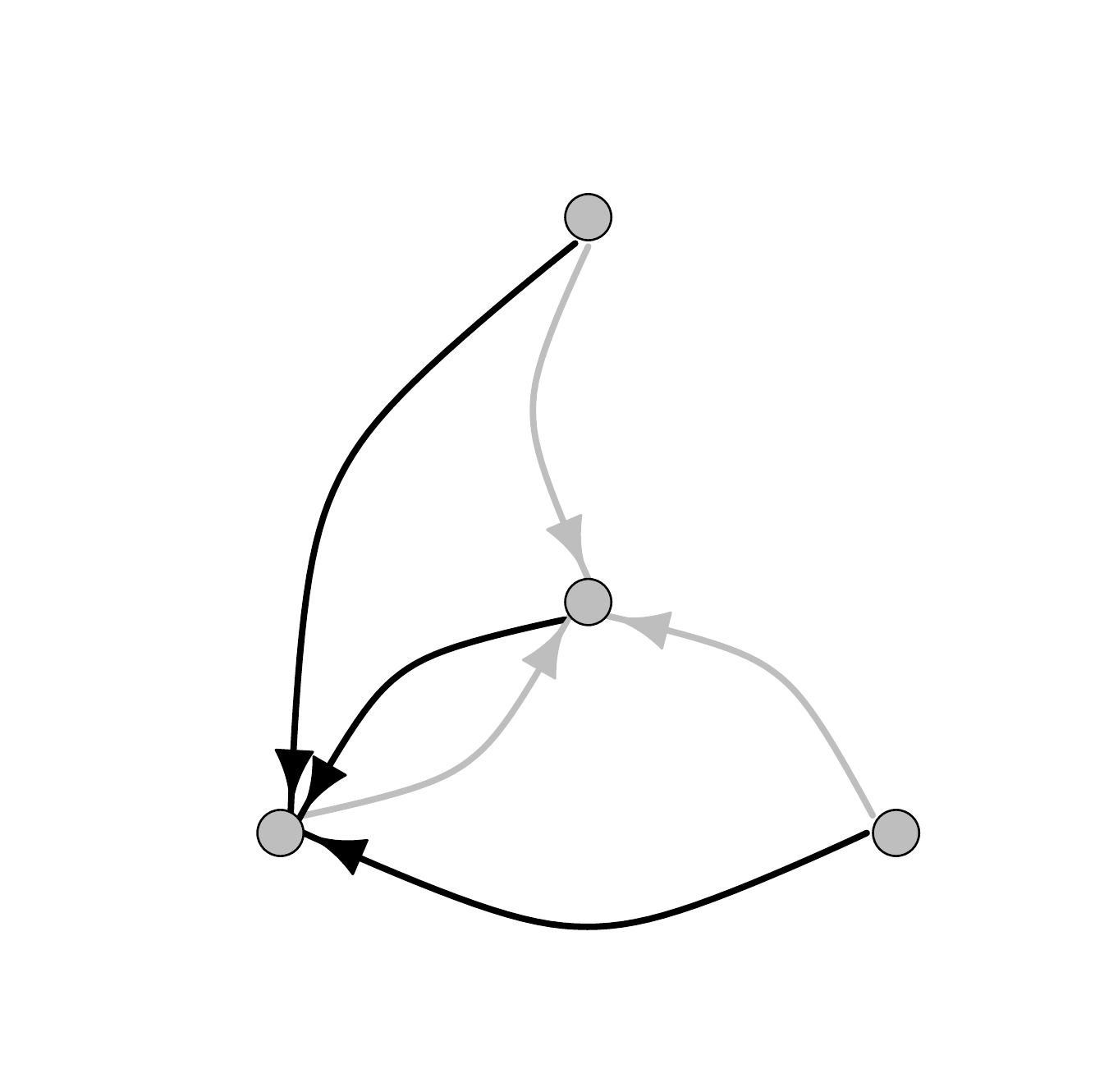}
    \caption{Case II/2 in the first step, $N=4$}\label{1p1q}
  \end{minipage}
\end{figure}

\setcounter{equation}{0}
First we give a precise mathematical description of our random graph model. 
For the convenience of the reader we also present figures. 
In Figures \ref{init}-\ref{1p1q} the initial and the first possible steps are shown when $N=4$ (the new vertices and edges are black and the old vertices and edges are grey). 

Let $N\geq 3$ be a fixed number.
We start at time $0$, and the initial graph is an $N$-star graph. 
Throughout the paper we call a graph $N$-star graph if it has $N$ vertices, 
one of them is called central vertex, the remaining $N-1$ vertices are called peripheral vertices, 
and it has $N-1$ directed edges.
The edges start from the $N-1$ peripheral vertices and their end point is the central vertex.
So the central vertex has in-degree $N-1$, and each of the $N-1$ peripheral vertices has out-degree $1$. 
The initial weight of the $N$-star is $1$, and the initial weights of its $(N-1)$-star sub-graphs are also $1$.
(The number of these $(N-1)$-star sub-graphs is $N-1$.)

The evolution of the graph is the following. 
At each step, $N$ vertices interact. 
Interaction (that is cooperation) means that we draw all edges from the peripheral vertices to the central vertex 
so that the vertices will form an $N$-star graph. 
We allow parallel edges.
When $N$ vertices interact, not only new edges are drawn, but the weights are also increased.
At the first interaction of $N$ vertices the newly created  $N$-star gets weight $1$, 
and its new $(N-1)$-star sub-graphs also get weight $1$.
If a sub-graph is not newly created, then its weight is increased by $1$.
When an existing $N$-star is activated again, then its weight and the weights of its $(N-1)$-star sub-graphs are increased by $1$.
So the weight of an $N$-star is the number of its activations.
We can see that the weight of an $(N-1)$-star is equal to the sum of the weights of the $N$-stars containing it.
The weights play crucial role in our model.
The higher the weight of a star the higher the chance that it will be activated again.

We have two options in every step of the evolution. 
{\tt Option I}:  with probability $p$, we add a new vertex, and it interacts with $N-1$ old vertices. 
{\tt Option II}: with probability $1-p$, we do not add any new vertex, but $N$ old vertices interact. 
Here $0<p\leq 1$ is fixed.

{\tt Option I.}
In this case, that is when a new vertex is born, we have again two possibilities: {\tt I/1} and {\tt I/2}.
\newline
{\tt I/1.}
The first possibility, which has probability $r$, is the following. (Here $0\leq r\leq 1$ is fixed.)
We choose one of the existing $(N-1)$-star sub-graphs according to the preferential attachment rule, 
and it will interact with the new vertex. 
Here the preferential attachment rule means that an $\left(N-1\right)$-star of weight $v_t$ is chosen with probability $v_t/\sum_h v_h$,
where $\sum_h v_h$ is the cumulated weight of the $(N-1)$-stars.
The interaction of the new vertex and the old $(N-1)$-star means that they establish a new $N$-star.
In this newly created $N$-star the center will be that vertex which was the center in the old $(N-1)$-star, 
the former $N-2$ peripheral vertices remain peripheral and the newly born vertex will be also peripheral.
New edges are drawn from the new and old peripheral vertices to the central one, and then the weights are increased by $1$.
More precisely, the just created $N$-star gets weight 1, among its $(N-1)$-star sub-graphs there are $(N-2)$ new ones, so each of them gets weight 1, 
finally the weight of the only old $(N-1)$-star sub-graph is increased by 1.
\newline
{\tt I/2.}
The second possibility has probability $1-r$.
Here we choose $N-1$ old vertices uniformly at random, and they will form an $N$-star graph with the new vertex, so that the new vertex will be the center. 
The edges are drawn from the peripheral vertices to the center.
Then the weights are increased.
As here the $N$-star graph and all of its $(N-1)$-star sub-graphs are new, so all of them get weight 1.

{\tt Option II.}
In this case, that is when we do not add any new vertex, we have two ways again: {\tt II/1} and {\tt II/2}.
\newline
{\tt II/1.}
The first way has probability $q$. (Here $0\leq q\leq 1$ is fixed.)
We choose one of the existing $N$-star sub-graphs by the preferential attachment rule, 
and draw all edges from its peripheral vertices to the center vertex. 
Then, as above, the weight of the $N$-star and the weights of its $(N-1)$-star sub-graphs are increased by $1$.
Here the preferential attachment rule means that an $N$-star of weight $v_t$ is chosen with probability $v_t/\sum_h v_h$,
where $\sum_h v_h$ is the cumulated weight of the $N$-stars.
\newline
{\tt II/2.}
The second way has probability $1-q$.  
We choose $N$ old vertices uniformly at random, and they establish an $N$-star graph.
Its center is chosen again uniformly at random out of the $N$ vertices. 
Then, as before, new edges are drawn from the peripheral vertices to the central one, and the weights of the $N$-star and its $(N-1)$-star sub-graphs are increased by $1$.

In this paper we show that this evolution leads to a scale-free graph. 


Throughout the paper $0<p\leq 1$,\;\; $0\leq r\leq 1$,\;\; $0\leq q\leq 1$ are fixed numbers. Let

\begin{equation*}
\alpha_{11}=pr,\;\;\;\;\alpha_{12}=(1-p)q,
\end{equation*}

\begin{equation*}
\alpha_1=\alpha_{11}+\alpha_{12},\;\;\;\;\alpha_2=pr\frac{N-2}{N-1}+(1-p)q,
\end{equation*}

\begin{equation}\label{parameters}
\beta_1=\frac{(1-p)(1-q)}{p},\;\;\;\;\beta_2=(N-1)\left[(1-r)+\frac{(1-p)(1-q)}{p}\right],
\end{equation}

\begin{equation*}
\alpha=\alpha_1+\alpha_2,\;\;\;\;\beta=\beta_1+\beta_2.
\end{equation*}

Let $V_n$ denote the number of vertices after $n$ steps. Let $Y \left( n,d_1,d_2 \right)$ denote the number of vertices with indegree $d_1$ and outdegree $d_2$ after the $n$th step. 
\begin{thm} \label{main}
Let $0<p<1$, $0<q<1$, $0<r<1$. Then for any fixed $d_1$ and $d_2$ with either $d_1=0$ and $1\leq d_2$ or $N-1\leq d_1$ and $d_2\geq 0$ we have
\begin{equation}
\dfrac{Y\left(n,d_1,d_2\right)}{V_n}\rightarrow y_{d_1,d_2}
\end{equation}
almost surely as $n\rightarrow\infty$, where $y_{d_1,d_2}$ are fixed non-negative numbers.


Let $d_2$ be fixed, then as $d_1\rightarrow\infty$
\begin{equation}
y_{d_1,d_2}\sim A(d_2)d_1^{-\left(1+\frac{\beta_2+1}{\alpha_1}\right)},
\end{equation}
where
\begin{equation}
A(d_2)=\dfrac{1-r}{\alpha_1}\dfrac{1}{d_2!}\dfrac{\varGamma\left(d_2+\dfrac{\beta_2}{\alpha_2}\right)}{\varGamma\left(\dfrac{\beta_2}{\alpha_2}\right)}\dfrac{\varGamma\left(1+\dfrac{\beta+1}{\alpha_1}\right)}{\varGamma\left(1+\dfrac{\beta_1}{\alpha_1}\right)}\dfrac{1}{\left(N-1\right)^{-\left(1+\frac{\beta_2+1}{\alpha_1}\right)}}.
\end{equation}

Let $d_1$ be fixed, then as $d_2\rightarrow\infty$
\begin{equation}
y_{d_1,d_2}\sim B(d_1)d_2^{-\left(1+\frac{\beta_1+1}{\alpha_2}\right)},
\end{equation}
where
\begin{equation}
B(d_1)=\dfrac{r}{\alpha_2}\dfrac{1}{\left(\dfrac{d_1}{N-1}\right)!}\dfrac{\varGamma\left(\dfrac{d_1}{N-1}+\dfrac{\beta_1}{\alpha_1}\right)}{\varGamma\left(\dfrac{\beta_1}{\alpha_1}\right)}\dfrac{\varGamma\left(1+\dfrac{\beta+1}{\alpha_2}\right)}{\varGamma\left(1+\dfrac{\beta_2}{\alpha_2}\right)}.
\end{equation}
Here $\varGamma$ denotes the Gamma function.
\end{thm}
\begin{pf}
The result is a consequence of Theorem \ref{thweight}.
\qed
\end{pf}
%
\section{Proofs and auxiliary results}  \label{proofs}
\setcounter{equation}{0}
%
\subsection{The evolution of the graph}
First we reformulate our model in order the simplify the proofs.
We shall see that, using the new parameters, our formulae will be symmetric, therefore we can shorten the proofs.
At the same time the new parametrization gives us a new viewpoint.
Our new description will be given in terms of undirected graphs without multiple edges.
Instead of the multiple edges we shall use weights of the vertices.
So we define for every vertex its central weight and its peripheral weight. 

The central weight of a vertex is $w_1$, if the vertex was $w_1$-times central vertex in interactions. 
The peripheral weight of a vertex is $w_2$, if the vertex was $w_2$-times peripheral vertex in interactions. 
It is easy to see that the central weight of a vertex is equal to $w_1= \dfrac{d_{1}}{N-1}$ and the peripheral weight of a vertex is equal to $w_2=d_{2}$, 
where $d_{1}$ denotes the in-degree of the vertex and $d_{2}$ denotes its out-degree. 
After the weights $w_1$ and $w_2$ are fixed, we delete all edges between any two given connected vertices and replace them by a single undirected edge.
Therefore this new edge will show that the two vertices cooperated at least once during the evolution.
In lemmas \ref{mainlemma}, \ref{xdw} and in theorems \ref{Xrelgyak}, \ref{thweight} we shall use this undirected graph.
We recall that the weight of an $N$-star is $w$, if the $N$-star took part in interactions $w$-times. 
Similarly, the weight of an $(N-1)$-star is $w$, if the $(N-1)$-star took part in interactions $w$-times. 
Recall that $V_n$ denotes the number of vertices after $n$ steps. 
Let $X \left( n,d,w_1,w_2 \right)$ denote the number of vertices of degree $d$, central weight $w_1$ and peripheral weight $w_2$ after the $n$th step. 
Furthermore, let $\FD_{n-1}$ denote the $\sigma$-algebra of observable events just after the $(n-1)$th step. 
We define $\binom{k}{l}=0$ if $l>k$.

\begin{lem}\label{mainlemma}
One has

\begin{equation*}
{\E} \{ X(n,d,w_1,w_2)|\FD_{n-1} \} =
X(n-1,d,w_1,w_2) \left[ 1-\left(\dfrac{w_1}{n}\alpha_1 + \dfrac{w_2}{n}\alpha_2 + \dfrac{p}{V_{n-1}}\beta \right) \right] +
\end{equation*}

\begin{equation*}
+X(n-1,d,w_1,w_2-1) \left[ pr \dfrac{(N-2)(w_2-1)}{(N-1)n} + \left( 1-p \right) \left( q\dfrac{w_2-1}{n} + \left( 1-q \right)\dfrac {d\binom{V_{n-1}-2}{N-2}} {\binom{V_{n-1}}{N} N}  \right) \right] +
\end{equation*}

\begin{equation*}
+X(n-1,d-1,w_1,w_2-1) \left[ p\left(1 - r\right)\dfrac{N-1}{V_{n-1}} + \left( 1-p \right) \left(1 - q\right) \dfrac {\left(V_{n-1}-d\right)\binom{V_{n-1}-2}{N-2}} {\binom{V_{n-1}}{N} N} \right] +
\end{equation*}

\begin{equation*}
+X(n-1,d,w_1-1,w_2) \left[ \left( 1-p \right) \left( q\dfrac{w_1-1}{n} + \left( 1-q \right)\dfrac {\binom{d}{N-1}} {\binom{V_{n-1}}{N} N} \right) \right] +
\end{equation*}

\begin{equation*}
+ X(n-1,d-1,w_1-1,w_2) \left[ pr \dfrac{ w_1-1 }{n} + \left( 1-p \right)\left( 1-q \right)\dfrac {\binom{d-1}{N-2} \left( V_{n-1}-d \right)} {\binom{V_{n-1}}{N}N} \right] +
\end{equation*}

\begin{equation*}
+ \sum_{m=2}^{N-2}X(n-1,d-m,w_1-1,w_2) \left[ \left( 1-p \right)\left( 1-q \right)\dfrac {\binom{d-m}{N-m-1} \binom{ V_{n-1}-d+m-1}{m} } {\binom{V_{n-1}}{N} N} \right]+
\end{equation*}

\begin{equation*} 
+X(n-1,d-\left(N-1\right),w_1-1,w_2) \left[ \left( 1-p \right) \left( 1-q \right)\dfrac{\binom{V_{n-1}-d+N-2}{N-1}}{\binom{V_{n-1}}{N} N} \right]+
\end{equation*}

\begin{equation}\label{ExpXndw}
+ pr\delta_{d,1}\delta_{w_1,0}\delta_{w_2,1}+p\left(1-r\right)\delta_{d,N-1}\delta_{w_1,1}\delta_{w_2,0}
\end{equation}

for either $w_1=0, 1 \leq w_2$ and $1 \leq d \leq w_2$ or $1 \leq w_1, 0 \leq w_2$ and $N-1 \leq d \leq w_1 (N - 1) + w_2$.
Here $\delta_{k,l}$ denotes the Dirac delta.
\end{lem}

\begin{pf}
Throughout the proof $w_1$ will denote the central weight and $w_2$ will denote the peripheral weight of a given vertex. 
The total weight of $(N-1)$-stars having a fixed common vertex of weights $w_1$ and $w_2$ is $w_1(N-1)+w_2(N-2)$. 
The total weight of $N$-stars after $(n-1)$ steps is $n$. The total weight of $(N-1)$-stars after $(n-1)$ steps is $n(N-1)$. 
 The probability that a given vertex is chosen, if we choose $(N-1)$ vertices uniformly is
\[
\frac{\binom{V_{n-1}-1}{N-2}}{\binom{V_{n-1}}{N-1}}=\frac{N-1}{V_{n-1}}.
\]
The probability that a given vertex is chosen, if we choose $N$ vertices uniformly is
\[
\frac{\binom{V_{n-1}-1}{N-1}}{\binom{V_{n-1}}{N}}=\frac{N}{V_{n-1}}.
\]  
So the probability that an old vertex of weights $w_1$ and $w_2$ takes part in the interaction at step $n$ is
\[
p\left(r\frac{w_1(N-1)+w_2(N-2)}{(N-1)n}+(1-r)\frac{N-1}{V_{n-1}}\right)+(1-p)\left(q\frac{w_1+w_2}{n}+(1-q)\frac{N}{V_{n-1}}\right)=
\]
\[
=\alpha_1\frac{w_1}{n}+\alpha_2\frac{w_2}{n}+\beta\frac{p}{V_{n-1}}.
\]
At each step when a new vertex is born we have two cases:
\begin{enumerate}
	\item with probability $pr$ a new vertex is born with central weight $0$, peripheral weight $1$ and degree $1$;
	\item with probability $p(1-r)$ a new vertex is born with central weight $1$, peripheral weight $0$ and degree $(N-1)$.
\end{enumerate}

Let us consider a fixed old vertex with degree $d$, central weight $w_1$ and peripheral weight $w_2$. For this vertex the probability that in the $n$th step
\begin{itemize}
\item neither its degree, nor its weights change is
\[
1-\left(\alpha_1\frac{w_1}{n}+\alpha_2\frac{w_2}{n}+\beta\frac{p}{V_{n-1}}\right);
\]
\item its degree does not change but its central weight is increased by $1$ is
\[
(1-p)\left(q\frac{w_1}{n}+(1-q)\frac{\binom{d}{N-1}}{N\binom{V_{n-1}}{N}}\right);
\]
\item its degree does not change but its peripheral weight is increased by $1$ is
\[
pr\frac{w_2(N-2)}{n(N-1)}+(1-p)q\frac{w_2}{n}+(1-p)(1-q)\frac{d\binom{V_{n-1}-2}{N-2}}{N\binom{V_{n-1}}{N}};
\]
\item its degree and its central weight are increased by $1$ is
\[
pr\frac{w_1}{n}+(1-p)(1-q)\frac{\binom{d}{N-2}\binom{V_{n-1}-d-1}{1}}{N\binom{V_{n-1}}{N}};
\]
\item its degree and its peripheral weight are increased by $1$ is
\[
p(1-r)\frac{N-1}{V_{n-1}}+(1-p)(1-q)\frac{\binom{V_{n-1}-d-1}{1}\binom{V_{n-1}-2}{N-2}}{N\binom{V_{n-1}}{N}};
\]
\item its degree is increased by $m$ $(1<m<N-1)$ and its central weight is increased by $1$ is
\[
(1-p)(1-q)\frac{\binom{V_{n-1}-d-1}{m}\binom{d}{N-1-m}}{N\binom{V_{n-1}}{N}};
\]
\item its degree is increased by $(N-1)$ and its central weight is increased by $1$ is
\[
(1-p)(1-q)\frac{\binom{V_{n-1}-d-1}{N-1}}{N\binom{V_{n-1}}{N}}.
\]
\end{itemize}
From these formulae, we obtain equation (\ref{ExpXndw}).
\qed
\end{pf}

\begin{thm}\label{Xrelgyak}
Let $0<p<1$, $0<q<1$, $0<r<1$.
Then for any fixed $w_1$, $w_2$ and $d$ with either $0=w_1$, $1\leq w_2$ and $1\leq d\leq w_2$ or $1 \leq w_1$, $0 \leq w_2$ and $N-1 \leq d \leq w_1\left(N-1 \right)+w_2$ we have
\begin{equation}   \label{X/Vtox}
\dfrac{X \left( n,d,w_1,w_2 \right)}{V_n} \rightarrow x_{d,w_1,w_2}
\end{equation}
almost surely as $n \rightarrow \infty $, where $x_{d,w_1,w_2}$ are fixed non-negative numbers.

Furthermore, the numbers $x_{d,w_1,w_2}$ satisfy the following recurrence relation
$$
x_{N-1,1,0} = \dfrac{1-r}{\alpha_1 + \beta +1} > 0, \quad\quad x_{d,1,0}=0, \, \text{ for} \quad d\ne N-1,
$$
$$
x_{1,0,1} = \dfrac{r}{\alpha_2 + \beta +1} > 0, \quad\quad x_{d,0,1}=0, \, \text{ for} \quad d\ne 1,
$$
\begin{equation} \label{rekurzio_x(d,w)}
x_{d,w_1,w_2} = \dfrac{1}{\alpha_1 w_1 + \alpha_2 w_2 + \beta +1} \times
\end{equation}
\begin{equation*}
\times \left[ \alpha_{11} \left( w_1-1 \right) x_{d-1,w_1-1,w_2} + \alpha_{12} \left(  w_1-1\right)x_{d,w_1-1,w_2} +\right.
\end{equation*}
\begin{equation*}
+\alpha_2\left(w_2-1\right)x_{d,w_1,w_2-1}+\left.\beta_1 x_{d-\left( N-1 \right),w_1-1,w_2} + \beta_2x_{d-1,w_1,w_2-1} \right]
\end{equation*}
for any $w_1$,$w_2$ and $d$. In the cases when $x_{d,w_1,w_2} = 0$ we have
$$
\dfrac{X \left( n,d,w_1,w_2 \right)}{V_n} = \o \left( n^{-a} \right),
$$
where $a$ is a positive number which may depend on $w_1$, $w_2$ and $d$.
\end{thm}

\begin{pf}
In the proof we shall use two major tools: a submartingale convergence theorem and mathematical induction.
For the reader's convenience we quote the submartingale convergence theorem in the Appendix (Theorem \ref{neveu}).
The proof will be divided into two main sections.
In the first section we shall introduce a submartingale and calculate its Doob-Meyer decomposition.
The second section will contain the mathematical induction.
As the index set is two-dimensional, we shall use multiple induction.
Therefore the second section of the proof will be divided again into subsections.

{\bf 1. The basic submartingale and its Doob-Meyer decomposition.}
Let
\begin{equation}\label{cnw1w2}
c(n,w_1,w_2)=\prod_{i=1}^n\left(1-\alpha_1\frac{w_1}{i}-\alpha_2\frac{w_2}{i}-\beta\frac{p}{V_{i-1}}\right)^{-1},\quad w_1+w_2\geq 1.
\end{equation}
We can see that
$c(n,w_1,w_2)$ is an $\FD_{n-1}$ measurable positive random variable. 
Using the Marcinkiewicz strong law of large numbers to the number of vertices,
we obtain that
\begin{equation}\label{V_nExp}
V_n=pn+\o\left(n^{1/2+\varepsilon}\right)
\end{equation}
almost surely, for any $\varepsilon>0$.

Now, using (\ref{V_nExp}) and the Taylor expansion for $\log(1+x)$, we have
\[
\log c(n,w_1,w_2)= - \sum_{i=1}^n \log\left(1-\alpha_1\frac{w_1}{i}-\alpha_2\frac{w_2}{i}-\beta\frac{1}{i+\o\left(i^{{1/2}+\varepsilon}\right)}\right)=
\]
\[
=\left(\alpha_1 w_1+\alpha_2 w_2+\beta\right)\sum_{i=1}^n\frac{1}{i}+\O(1),
\]
where the error term is convergent as $n\rightarrow\infty$. So
\begin{equation}\label{cnw1w2asmp}
c(n,w_1,w_2)\sim a_{w_1,w_2}n^{\alpha_1 w_1+\alpha_2 w_2+\beta}
\end{equation}
almost surely, as $n\rightarrow\infty$. Here $a_{w_1,w_2}$ is a positive random variable.

Let
\begin{equation*}
Z(n,d,w_1,w_2)=c(n,w_1,w_2)X(n,d,w_1,w_2),
\end{equation*}
where $1\leq d\leq w_1(N-1)+w_2,1\leq w_1+w_2$.
In formula \eqref{ExpXndw} all terms are non-negative.
So multiplying both sides of \eqref{ExpXndw} by $c(n,w_1,w_2)$, we see that
 $\{Z(n,d,w_1,w_2),\FD_n,n=1,2,\dots\}$ is a non-negative submartingale for any fixed $1\leq d\leq w_1(N-1)+w_2$ and $1\leq w_1+w_2$. 
 By the Doob-Meyer decomposition of $Z(n,d,w_1,w_2)$, we have
\begin{equation*}
Z(n,d,w_1,w_2)=M(n,d,w_1,w_2)+A(n,d,w_1,w_2),
\end{equation*}
where $M(n,d,w_1,w_2)$ is a martingale, and $A(n,d,w_1,w_2)$ is a predictable increasing process, and their general forms are the following

\begin{equation}\label{Mgform}
M(n,d,w_1,w_2)=\sum_{i=1}^n\left[Z(i,d,w_1,w_2)-\E(Z(i,d,w_1,w_2)|\FD_{i-1})\right],
\end{equation}

\begin{equation}\label{Agform}
A(n,d,w_1,w_2)=\E Z(1,d,w_1,w_2)+\sum_{i=2}^n\left[\E(Z(i,d,w_1,w_2)|\FD_{i-1})-Z(i-1,d,w_1,w_2)\right].
\end{equation}
Here $\FD_0$ denotes the trivial $\sigma$-algebra. 
From Lemma \ref{mainlemma} and (\ref{Agform}), we obtain

\begin{equation*}
A(n,d,w_1,w_2)=\E Z(1,d,w_1,w_2)+
\end{equation*}

\begin{equation*}
+\sum_{i=2}^n c(i,w_1,w_2)\Bigg\{ X(i-1,d,w_1,w_2-1) \left[ pr \dfrac{(N-2)(w_2-1)}{(N-1)i} +\right. \Bigg.
\end{equation*}

\begin{equation*}
+\left.\left( 1-p \right) \left( q\dfrac{w_2-1}{i} + \left( 1-q \right)\dfrac {d\binom{V_{i-1}-2}{N-2}} {\binom{V_{i-1}}{N} N}  \right) \right] +
\end{equation*}

\begin{equation*}
+X(i-1,d-1,w_1,w_2-1) \left( p\left(1 - r\right)\dfrac{N-1}{V_{i-1}} + \left( 1-p \right) \left(1 - q\right) \dfrac {\left(V_{i-1}-d\right)\binom{V_{i-1}-2}{N-2}} {\binom{V_{i-1}}{N} N} \right)+
\end{equation*}

\begin{equation}\label{Aform}
+X(i-1,d,w_1-1,w_2) \left(\left( 1-p \right)  q\dfrac{w_1-1}{i}\right) + X(i-1,d-1,w_1-1,w_2) \left(pr \dfrac{ w_1-1 }{i}\right)+
\end{equation}

\begin{equation*}
+ \sum_{m=0}^{N-1}X(i-1,d-m,w_1-1,w_2) \left( \left( 1-p \right)\left( 1-q \right)\dfrac {\binom{d-m}{N-m-1} \binom{ V_{i-1}-d+m-1}{m} } {\binom{V_{i-1}}{N} N} \right)+
\end{equation*}

\begin{equation*}
+ \Bigg. pr\delta_{d,1}\delta_{w_1,0}\delta_{w_2,1}+p\left(1-r\right)\delta_{d,N-1}\delta_{w_1,1}\delta_{w_2,0}\Bigg\}.
\end{equation*}

In the following we give an upper bound for $B(n,d,w_1,w_2)$, where $B(n,d,w_1,w_2)$ denotes the sum of the conditional variances of $Z(n,d,w_1,w_2)$.

\begin{equation}\label{Bform}
B(n,d,w_1,w_2)=\sum_{i=2}^n\E\left\{(Z(i,d,w_1,w_2)-\E(Z(i,d,w_1,w_2)|\FD_{i-1}))^2|\FD_{i-1}\right\}=
\end{equation}

\begin{equation*}
=\sum_{i=2}^n c(i,w_1,w_2)^2\E\left\{(X(i,d,w_1,w_2)-\E(X(i,d,w_1,w_2)|\FD_{i-1}))^2|\FD_{i-1}\right\}\leq
\end{equation*}

\begin{equation*}
\leq\sum_{i=2}^n c(i,w_1,w_2)^2\E\left\{(X(i,d,w_1,w_2)-X(i-1,d,w_1,w_2))^2|\FD_{i-1}\right\}\leq
\end{equation*}

\begin{equation*}
\leq N^2\sum_{i=2}^n c(i,w_1,w_2)^2=\O \left( n^{2(\alpha_1w_1+\alpha_2w_2+\beta)+1}\right).
\end{equation*}
Here first we used that $c(i,w_1,w_2)$ is $\FD_{i-1}$ measurable, then the fact that at each step $N$ vertices interact, finally we applied (\ref{cnw1w2asmp}).

As $M(n,d,w_1,w_2)$ is a martingale, therefore $M^2(n,d,w_1,w_2)$ is a submartingale according to Jensen's inequality. 
Applying the Doob-Meyer decomposition for $M^2(n,d,w_1,w_2)$, we obtain

\begin{equation}\label{M2form}
M^2(n,d,w_1,w_2)=Y(n,d,w_1,w_2)+B(n,d,w_1,w_2),
\end{equation}
where $Y(n,d,w_1,w_2)$ is a martingale and the predictable increasing process $B(n,d,w_1,w_2)$ is the same as the one in (\ref{Bform}).

{\bf 2. The mathematical induction.}
First we consider the particular case $w_1=1$, $w_2=0$, then the case $w_1=0$, $w_2=1$ and also the case $w_1=1$, $w_2=1$. 
Then we use induction along the boundary of the domain, i.e. when $w_1=0$ or $w_2=0$.
Finally, we use induction in the interior of the domain.

{\bf Step 2/a.}
Let $w_1=1$ and $w_2=0$. 
A vertex with these weights exists if and only if it was center once and it was not even once peripheral. 
In this case its degree has to be equal to $N-1$. 
Using (\ref{Aform}) and (\ref{cnw1w2asmp}),

\begin{equation*}
A(n,N-1,1,0)\sim p(1-r)\sum_{i=2}^n c(i,1,0) \sim p(1-r)\sum_{i=2}^n a_{1,0} i^{\alpha_1+\beta} \sim
\end{equation*}

\begin{equation}\label{A10}
\sim p(1-r)a_{1,0}\frac{n^{\alpha_1+\beta+1}}{\alpha_1+\beta+1}\rightarrow\infty,
\end{equation}
almost surely as $n\rightarrow\infty$.
Using (\ref{Bform}),
\begin{equation*}
B(n,N-1,1,0)=\O(n^{2(\alpha_1+\beta)+1}),
\end{equation*}
so
\begin{equation*}
B(n,N-1,1,0)^{\frac{1}{2}}\log B(n,N-1,1,0)=\O(A(n,N-1,1,0)).
\end{equation*}
Now, applying Lemma \ref{neveu}, we get
\begin{equation}\label{Z10}
Z(n,N-1,1,0)\sim A(n,N-1,1,0)
\end{equation}
almost surely on the event $\{A(n,N-1,1,0)\rightarrow\infty\}$ as $n\rightarrow\infty$.
We have, by using (\ref{A10}), (\ref{V_nExp}), (\ref{cnw1w2asmp}) and (\ref{Z10}), that
\begin{equation*}
\frac{X(n,N-1,1,0)}{V_n}=\frac{Z(n,N-1,1,0)}{c(n,1,0)V_n}\sim\frac{A(n,N-1,1,0)}{c(n,1,0)V_n}\sim
\end{equation*}
\begin{equation}\label{xN110}
\sim\frac{p(1-r)a_{1,0}\frac{n^{\alpha_1+\beta+1}}{\alpha_1+\beta+1}}{a_{1,0}n^{\alpha_1+\beta}pn}=\frac{1-r}{\alpha_1+\beta+1}=x_{N-1,1,0}> 0.
\end{equation}
Otherwise, when $d\neq N-1$, then $X(n,d,1,0)\equiv 0$. So (\ref{rekurzio_x(d,w)}) is true for $w_1=1$ and $w_2=0$.

{\bf Step 2/b.}
Now let $w_1=0$ and $w_2=1$. 
In this case the degree of the vertex has to be $1$. From (\ref{Aform}) 
\begin{equation}
A(n,1,0,1)\sim pr\sum_{i=2}^n c(i,0,1) \sim pr\sum_{i=2}^n a_{0,1} i^{\alpha_2+\beta} \sim
\end{equation}

\begin{equation}\label{A01}
\sim pra_{0,1}\frac{n^{\alpha_2+\beta+1}}{\alpha_2+\beta+1}\rightarrow\infty,
\end{equation}
almost surely as $n\rightarrow\infty$.
Using (\ref{Bform}), we have
\begin{equation*}
B(n,1,0,1)=\O(n^{2(\alpha_1+\beta)+1}),
\end{equation*}
so
\begin{equation*}
B(n,1,0,1)^{\frac{1}{2}}\log B(n,1,0,1)=\O(A(n,1,0,1)).
\end{equation*}
As before, Lemma \ref{neveu}, we obtain 
\begin{equation}\label{Z01}
Z(n,1,0,1)\sim A(n,1,0,1)
\end{equation}
almost surely on the event $\{A(n,1,0,1)\rightarrow\infty\}$ as $n\rightarrow\infty$.
So, by using (\ref{A01}), (\ref{V_nExp}), (\ref{cnw1w2asmp}) and (\ref{Z01}), we have
\begin{equation*}
\frac{X(n,1,0,1)}{V_n}=\frac{Z(n,1,0,1)}{c(n,0,1)V_n}\sim\frac{A(n,1,0,1)}{c(n,0,1)V_n}\sim
\end{equation*}
\begin{equation}\label{x101}
\sim\frac{pra_{0,1}\frac{n^{\alpha_2+\beta+1}}{\alpha_2+\beta+1}}{a_{0,1}n^{\alpha_2+\beta}pn}=\frac{r}{\alpha_2+\beta+1}=x_{1,0,1}> 0.
\end{equation}
Otherwise, when  $d\neq 1$, then $X(n,d,0,1)\equiv 0$. So (\ref{rekurzio_x(d,w)}) is true for $w_1=0$ and $w_2=1$, too.

{\bf Step 2/c.}
Now consider the case of $w_1=w_2=1$. 
Vertices with these weights exist only with degree $N-1$ or $N$. 
First we deal with the case, when the degree is $N-1$. 
Now by (\ref{Aform}), (\ref{cnw1w2asmp}), (\ref{xN110}) and (\ref{x101}),

\begin{equation*}
A(n,N-1,1,1)\sim
\end{equation*}

\begin{equation*}
\sim\sum_{i=2}^n a_{1,1}i^{\alpha_1+\alpha_2+\beta}\left(x_{1,0,1}(1-p)(1-q)\frac{N-1}{p(i-1)}+x_{N-1,1,0}(1-p)(1-q)\frac{(N-1)^2}{p(i-1)}\right)\sim
\end{equation*}

\begin{equation}\label{AN111}
\sim\sum_{i=2}^ni^{\alpha_1+\alpha_2+\beta-1}a_{1,1}\frac{(1-p)(1-q)(N-1)}{p}\left(x_{1,0,1}+x_{N-1,1,0}(N-1)\right)\sim
\end{equation}

\begin{equation*}
\sim a_{1,1}\frac{(1-p)(1-q)(N-1)}{p}\left(x_{1,0,1}+(N-1)x_{N-1,1,0}\right)\frac{n^{\alpha_1+\alpha_2+\beta}}{\alpha_1+\alpha_2+\beta}.
\end{equation*}

By (\ref{Bform})
\begin{equation*}
B(n,N-1,1,1)=\O(n^{2(\alpha_1+\alpha_2+\beta)+1}).
\end{equation*}
So in this case
\begin{equation*}
B(n,N-1,1,1)^{\frac{1}{2}}\log B(n,N-1,1,1)\neq \O(A(n,N-1,1,1)).
\end{equation*}
Using Lemma \ref{neveu}, we can see that for any $\varepsilon>0$

\begin{equation*}
M(n,N-1,1,1)=\o(B(n,N-1,1,1)^{\frac{1}{2}}\log B(n,N-1,1,1))=\o(n^{\alpha_1+\alpha_2+\beta+\frac{1}{2}+\varepsilon})
\end{equation*}
almost surely on the event $\{B(n,N-1,1,1)\rightarrow\infty\}$ as $n\rightarrow\infty$. 
Furthermore $M(n,N-1,1,1)$ is convergent on the event $\{B(\infty,N-1,1,1)<\infty\}$ as $n\rightarrow\infty$, 
so $M(n,N-1,1,1)=\o(n^{\alpha_1+\alpha_2+\beta+\frac{1}{2}+\varepsilon})$ almost surely. 
Therefore, from (\ref{AN111}), (\ref{cnw1w2asmp}) and (\ref{V_nExp}), we obtain

\begin{equation*}
\frac{X(n,N-1,1,1)}{V_n}=\frac{Z(n,N-1,1,1)}{c(n,1,1)V_n}=\frac{M(n,N-1,1,1)+A(n,N-1,1,1)}{c(n,1,1)V_n}\leq
\end{equation*}

\begin{equation}
\leq\frac{Cn^{\alpha_1+\alpha_2+\beta+\frac{1}{2}+\varepsilon}}{n^{\alpha_1+\alpha_2+\beta}n}\leq C\frac{1}{n^a}\rightarrow 0,
\end{equation}
 as $n\rightarrow\infty$ with $\frac{1}{4}<a<\frac{1}{2}$. 

Consider the second case, that is when the degree of vertices is equal to $N$. 
From (\ref{Aform})
 
\begin{equation*}
A(n,N,1,1)\sim
\end{equation*}

\begin{equation*}
\sim\sum_{i=2}^n a_{1,1} i^{\alpha_1+\alpha_2+\beta}\left[(N-1)\left(p(1-r)+(1-p)(1-q)\right)x_{N-1,1,0}+(1-p)(1-q)x_{1,0,1}\right]\sim
\end{equation*}

\begin{equation}\label{AN11}
\sim a_{1,1}\left[(N-1)\left(p(1-r)+(1-p)(1-q)\right)x_{N-1,1,0}+\right.
\end{equation}

\begin{equation*}
\left.+(1-p)(1-q)x_{1,0,1}\right]\frac{n^{\alpha_1+\alpha_2+\beta+1}}{\alpha_1+\alpha_2+\beta+1}\rightarrow\infty
\end{equation*}
almost surely as $n\rightarrow\infty$.
Using (\ref{Bform}),
\begin{equation*}
B(n,N,1,1)=\O(n^{2(\alpha_1+\alpha_2+\beta)+1}),
\end{equation*}
so
\begin{equation*}
B(n,N,1,1)^{\frac{1}{2}}\log B(n,N,1,1)=\O(A(n,N,1,1)).
\end{equation*}
Similarly as before, using Lemma \ref{neveu}, we obtain 
\begin{equation}\label{ZN11}
Z(n,N,1,1)\sim A(n,N,1,1)
\end{equation}
almost surely on the event $\{A(n,N,1,1)\rightarrow\infty\}$ as $n\rightarrow\infty$.
So from (\ref{AN11}), (\ref{V_nExp}), (\ref{cnw1w2asmp}) and (\ref{ZN11}) we obtain
\begin{equation*}
\frac{X(n,N,1,1)}{V_n}=\frac{Z(n,N,1,1)}{c(n,1,1)V_n}\sim\frac{A(n,N,1,1)}{c(n,1,1)V_n}\sim
\end{equation*}

\begin{equation*}
\sim\frac{(N-1)\left(p(1-r)+(1-p)(1-q)\right)x_{N-1,1,0}+(1-p)(1-q)x_{1,0,1}}{(\alpha_1+\alpha_2+\beta+1)p}=
\end{equation*}

\begin{equation}
=\frac{1}{\alpha_1+\alpha_2+\beta+1}\left(\beta_1x_{1,0,1}+\beta_2x_{N-1,1,0}\right)=x_{N,1,1}> 0.
\end{equation}
Otherwise, if $d\neq N-1$ or $d\neq N$, then $X(n,d,1,1)\equiv 0$. So (\ref{rekurzio_x(d,w)}) is true for $w_1=1$ and $w_2=1$, too.

{\bf Step 2/d.}
Now we study the case of $w_1=k$ and $w_2=0$, $k>1$. 
These vertices were always central in interactions, and they never were peripheral. 
In this case $N-1\leq d\leq k(N-1)$, where $d$ denotes the degree of vertices. 
Suppose that the statement is true for all central weights less than $k$, for zero peripheral weight and for all possible degrees. 
Assume that at least one of the coefficients $x_{d-1,k-1,0}, x_{d,k-1,0}, x_{d-(N-1),k-1,0}$ is positive. As earlier, from (\ref{Aform})

\begin{equation*}
A(n,d,k,0)\sim\sum_{i=2}^n a_{k,0}i^{k\alpha_1+\beta}\left[ x_{d-1,k-1,0}(k-1)p^2r+x_{d,k-1,0}(k-1)p(1-p)q\right.+
\end{equation*}

\begin{equation*}
+\left. x_{d-(N-1),k-1,0}(1-p)(1-q)\right]\sim
\end{equation*}

\begin{equation}\label{Adk0}
\sim a_{k,0}\frac{n^{k\alpha_1+\beta+1}}{k\alpha_1+\beta+1}\left[ x_{d-1,k-1,0}(k-1)p^2r+x_{d,k-1,0}(k-1)p(1-p)q+\right.
\end{equation}

\begin{equation*}
\left. +x_{d-(N-1),k-1,0}(1-p)(1-q)\right]\rightarrow\infty,
\end{equation*}
almost surely as $n\rightarrow\infty$.

Using (\ref{Bform}), we have

\begin{equation}\label{Bndk0}
B(n,d,k,0)=\O(n^{2(k\alpha_1+\beta)+1}),
\end{equation}
so
\begin{equation}\label{Bndk00}
B(n,d,k,0)^\frac{1}{2}\log B(n,d,k,0)=\O(A(n,d,k,0)).
\end{equation}
Similarly as above, we have

\begin{equation}\label{Zdk0}
Z(n,d,k,0)\sim A(n,d,k,0)
\end{equation}
almost surely on the event $\{A(n,d,k,0)\rightarrow\infty\}$ as $n\rightarrow\infty$. 
Therefore, by using (\ref{Adk0}), (\ref{V_nExp}), (\ref{cnw1w2asmp}) and (\ref{Zdk0}), we have
\begin{equation*}
\frac{X(n,d,k,0)}{V_n}=\frac{Z(n,d,k,0)}{c(n,k,0)V_n}\sim \frac{A(n,d,k,0)}{c(n,k,0)V_n}\sim
\end{equation*}

\begin{equation}
\sim\frac{1}{k\alpha_1+\beta+1}\left[ (k-1) \left(\alpha_{11}x_{d-1,k-1,0}+\alpha_{12}x_{d,k-1,0}\right)+\beta_1x_{d-(N-1),k-1,0}\right]=x_{d,k,0}>0,
\end{equation}
if at least one of the coefficients $x_{d-1,k-1,0}$, $x_{d,k-1,0}$ and $x_{d-(N-1),k-1,0}$ is positive. Otherwise, the limit is zero. One can see for any fixed $k$ the limit of $\frac{X(n,d,k,0)}{V_n}$ is positive if $d\geq N-1$ and $d$ is 'close to' $N-1$.

However, if $d$ is 'close to' $k(N-1)$, the limit of $\frac{X(n,d,k,0)}{V_n}$ can be either zero or positive. If $d< N-1$ or $d> k(N-1)$, the limit is zero, because in these cases $X(n,d,k,0)\equiv 0$. 

Now, consider the case, when the coefficients $x_{d-1,k-1,0}, x_{d,k-1,0}$, and $x_{d-(N-1),k-1,0}$ are equal to zero. 
From (\ref{Aform}), (\ref{V_nExp}), (\ref{cnw1w2asmp}) and using the induction hypothesis, we obtain

\begin{multline*}
A(n,d,k,0)\sim \sum_{i=2}^n a_{k,0}i^{k\alpha_1+\beta}\left(\O\left(\dfrac{1}{i^a}\right) + \sum_{m=2}^{N-2}X(i-1,d-m,k-1,0) \times\right. \\\left. \times \left( \left( 1-p \right)\left( 1-q \right)\dfrac {\binom{d-m}{N-m-1} (N-1)! } {m!}\dfrac{1}{(pi)^{N-m}} \right)\right)\leq
\end{multline*}
\[
\leq C_1\sum_{i=2}^n a_{k,0}i^{k\alpha_1+\beta-a}+C_2 \sum_{i=2}^n\sum_{m=2}^{N-2}a_{k,0}i^{k\alpha_1+\beta}x_{d-m,k-1,0}\dfrac{1}{i^{N-m-1}}\leq
\]
\begin{equation}\label{Adk00}
\leq a_{k,0}C_1\dfrac{n^{k\alpha_1+\beta+1-a}}{k\alpha_1+\beta+1-a}+a_{k,0}C_2\dfrac{n^{k\alpha_1+\beta}}{k\alpha_1+\beta}=\O\left( n^{k\alpha_1+\beta+1-a}\right),
\end{equation}
where $C_1$ and $C_2$ appropriate constants. 
So we can not apply (\ref{Bndk00}). 
But, applying (\ref{V_nExp}), (\ref{cnw1w2asmp}) and Lemma \ref{neveu}, we obtain
\[
\dfrac{X(n,d,k,0)}{V_n}=\dfrac{Z(n,d,k,0)}{c(n,k,0)V_n}=\dfrac{M(n,d,k,0)+A(n,d,k,0)}{c(n,k,0)V_n}=
\]
\begin{equation}
=\dfrac{\O\left(n^{k\alpha_1+\beta+1-a}\right)}{n^{k\alpha_1+\beta}pn}=\O\left(n^{-a}\right)\rightarrow 0 = x_{d,k,0}
\end{equation}
almost surely as $n\rightarrow\infty$.

{\bf Step 2/e.}
Consider the case of $w_1=0$ and $w_2=l$, $l>1$. 
These vertices were always peripheral in the interactions, they never were center. 
In this case $1\leq d\leq l$, where $d$ denotes the degree of vertices. 
Suppose that the statement is true for all peripheral weights less than $l$, for central weight zero, and for all possible degrees. 
One can see that at least one of the coefficients $x_{d-1,0,l-1}, x_{d,0,l-1}$ is positive. 
As before, from (\ref{Aform}) we obtain that

\begin{equation*}
A(n,d,0,l)\sim\sum_{i=2}^na_{0,l}i^{l\alpha_2+\beta}\left[x_{d,0,l-1}p(l-1)\left(\frac{N-2}{N-1}pr+(1-p)q\right)+ \right.
\end{equation*}

\begin{equation*}
\left.+x_{d-1,0,l-1}(N-1)\left(p(1-r)+(1-p)(1-q)\right)\frac{}{}\right]\sim
\end{equation*}

\begin{equation}\label{Ad0l}
\sim a_{0,l}\frac{n^{l\alpha_2+\beta+1}}{l\alpha_2+\beta+1}\left[x_{d,0,l-1}p(l-1)\left(\frac{N-2}{N-1}pr+(1-p)q\right)+ \right.
\end{equation}

\begin{equation*}
\left.+x_{d-1,0,l-1}(N-1)\left(p(1-r)+(1-p)(1-q)\right)\frac{}{}\right]\rightarrow\infty,
\end{equation*}
almost surely as $n\rightarrow\infty$. Applying (\ref{Bform}), we have

\begin{equation*}
B(n,d,0,l)=\O(n^{2(l\alpha_2+\beta)+1}),
\end{equation*}
therefore

\begin{equation*}
B(n,d,0,l)^{\frac{1}{2}}\log B(n,d,0,l)=\O(A(n,d,0,l)).
\end{equation*}
Using again Lemma \ref{neveu}, we see that

\begin{equation}
Z(n,d,0,l)\sim A(n,d,0,l)
\end{equation}
almost surely on the event $\{A(n,d,0,l)\rightarrow\infty\}$ as $n\rightarrow\infty$. From the above formulae, we obtain

\begin{equation*}
\frac{X(n,d,0,l)}{V_n}=\frac{Z(n,d,0,l)}{c(n,0,l)V_n}\sim \frac{A(n,d,0,l)}{c(n,0,l)V_n}\sim
\end{equation*}

\begin{equation}
\sim\frac{1}{l\alpha_2+\beta+1}\left[ \alpha_2(l-1)x_{d,0,l-1}+\beta_2x_{d-1,0,l-1}\right]=x_{d,0,l}>0,
\end{equation}
because at least one of the coefficients $x_{d-1,0,l-1}, x_{d,0,l-1}$ is positive. If $d<1$ or $d>l$, then $X(n,d,0,l)\equiv 0$.

{\bf Step 2/f.}
Finally, the last part of our proof, is the case when $w_1=k>0$ and $w_2=l>0$. 
These vertices were center $k$ times and were peripheral $l$ times in the interactions. 
Here $N-1\leq d\leq k(N-1)+l$. 
Suppose that the statement is true for any vertex if either its central weight is less than $k$ and its peripheral weight is not greater than $l$ 
or its central weight is not greater than $k$ and its peripheral weight is less than $l$. 
Assume that at least one of the coefficients $x_{d-1,k-1,l}, x_{d,k,l-1}, x_{d-1,k,l-1}, x_{d,k-1,l}, x_{d-(N-1),k-1,l}$ is positive. 
From (\ref{Aform}) and using the induction hypothesis, we get

\begin{equation*}
A(n,d,k,l)\sim\sum_{i=2}^n a_{k,l}i^{k\alpha_1+l\alpha_2+\beta}\left[x_{d-1,k-1,l}(k-1)p^2r+\right.
\end{equation*}

\begin{equation*}
+x_{d-1,k,l-1}(N-1)\left(p(1-r)+(1-p)(1-q)\right)+x_{d,k-1,l}(k-1)p(1-p)q+
\end{equation*}

\begin{equation*}
\left.+ x_{d,k,l-1}(l-1)p\left(pr\frac{N-2}{N-1}+(1-p)q\right)+x_{d-(N-1),k-1,l}(1-p)(1-q)\frac{}{}\right]\sim
\end{equation*}

\begin{equation*}
\sim a_{k,l}\frac{n^{k\alpha_1+l\alpha_2+\beta +1}}{k\alpha_1+l\alpha_2+\beta +1}\left[x_{d-1,k-1,l}(k-1)p^2r\right.+
\end{equation*}

\begin{equation}\label{Adkl}
+x_{d-1,k,l-1}(N-1)\left(p(1-r)+(1-p)(1-q)\right)+x_{d,k-1,l}(k-1)p(1-p)q+
\end{equation}

\begin{equation*}
\left.+x_{d,k,l-1}(l-1)p\left(pr\frac{N-2}{N-1}+(1-p)q\right)+ x_{d-(N-1),k-1,l}(1-p)(1-q)\right]\rightarrow\infty
\end{equation*}
almost surely as $n\rightarrow\infty$. Applying (\ref{Bform}), we obtain

\begin{equation}\label{Bndkl}
B(n,d,k,l)=\O(n^{2(k\alpha_1+l\alpha_2+\beta) +1}),
\end{equation}
so

\begin{equation*}
B(n,d,k,l)^{\frac{1}{2}}\log B(n,d,k,l) =\O(A(n,d,k,l)).
\end{equation*}

Using Lemma \ref{neveu}, we get

\begin{equation}\label{Zdkl}
Z(n,d,k,l)\sim A(n,d,k,l)
\end{equation}
almost surely on the event $\{A(n,d,k,l)\rightarrow\infty\}$ as $n\rightarrow\infty$. Therefore, from (\ref{Adkl}), (\ref{V_nExp}), (\ref{cnw1w2asmp}) and (\ref{Zdkl}),

\begin{equation*}
\frac{X(n,d,k,l)}{V_n}=\frac{Z(n,d,k,l)}{c(n,k,l)V_n}\sim \frac{A(n,d,k,l)}{c(n,k,l)V_n}\sim
\end{equation*}

\begin{equation*}
\sim\frac{1}{k\alpha_1+l\alpha_2+\beta+1}\left[\alpha_{11}(k-1)x_{d-1,k-1,l}+\alpha_{12}(k-1)x_{d,k-1,l}+\alpha_2(l-1)x_{d,k,l-1}+\right.
\end{equation*}

\begin{equation}
\left.+\beta_1x_{d-(N-1),k-1,l}+\beta_2x_{d-1,k,l-1}\right]=x_{d,k,l}>0,
\end{equation}
because of the assumptions. 
Otherwise, if $d<N-1$ or $d>k(N-1)+l$, then $X(n,d,k,l)\equiv 0$. 
So in this case (\ref{rekurzio_x(d,w)}) is true for all $k,l\in\mathbb{N}$.

The remaining case is when the coefficients $x_{d-1,k-1,l}$, $x_{d,k,l-1}$, $x_{d-1,k,l-1}$, $x_{d,k-1,l}$, and $x_{d-(N-1),k-1,l}$ are all equal to zero.
From (\ref{Aform}) and using the induction hypothesis, we obtain
\begin{multline*}
A(n,d,k,l)\sim \sum_{i=2}^n a_{k,l}i^{k\alpha_1+l\alpha_2+\beta}\left(\O\left(\dfrac{1}{i^a}\right) + \sum_{m=2}^{N-2}X(i-1,d-m,k-1,l) \times\right. \\\left. \times \left( \left( 1-p \right)\left( 1-q \right)\dfrac {\binom{d-m}{N-m-1} (N-1)! } {m!}\dfrac{1}{(pi)^{N-m}} \right)\right)\sim
\end{multline*}
\[
\sim C_1\sum_{i=2}^n a_{k,l}i^{k\alpha_1+l\alpha_2+\beta-a}+C_2 \sum_{i=2}^n\sum_{m=2}^{N-2}a_{k,l}i^{k\alpha_1+l\alpha_2+\beta}x_{d-m,k-1,l}\dfrac{1}{i^{N-m-1}}\leq
\]
\[
\leq a_{k,l}C_1\dfrac{n^{k\alpha_1+l\alpha_2+\beta+1-a}}{k\alpha_1+l\alpha_2+\beta+1-a}+a_{k,l}C_2\dfrac{n^{k\alpha_1+l\alpha_2+\beta}}{k\alpha_1+l\alpha_2+\beta}=
\]
\begin{equation}\label{Adkl0}
=\O\left( n^{k\alpha_1+l\alpha_2+\beta+1-a}\right),
\end{equation}
where $C_1$ and $C_2$ are appropriate constants. 
Using (\ref{V_nExp}), (\ref{cnw1w2asmp}) and Lemma \ref{neveu}, we obtain
\[
\dfrac{X(n,d,k,l)}{V_n}=\dfrac{Z(n,d,k,l)}{c(n,k,l)V_n}=\dfrac{M(n,d,k,l)+A(n,d,k,l)}{c(n,k,l)V_n}=
\]
\begin{equation}
=\dfrac{\O\left(n^{k\alpha_1+l\alpha_2+\beta+1-a}\right)}{n^{k\alpha_1+l\alpha_2+\beta}pn}=\O\left(n^{-a}\right)\rightarrow 0 = x_{d,k,l}
\end{equation}
almost surely as $n\rightarrow\infty$.

\qed
\end{pf}

\subsection{Scale-free property of the weights}

\begin{lem}  \label{xdw}
Let $p>0$ and define
\begin{equation*}
x_{w_1,w_2} = x_{1,w_1,w_2} + x_{2,w_1,w_2} + \dots + x_{\left(N-1\right)w_1+w_2,w_1,w_2}
\end{equation*}
for $w_1\ge 0$, $w_2\ge 0$ and $w_1+w_2\ge 1$.
Then $x_{w_1,w_2}$ ($1\leq w_1+w_2$) are positive numbers satisfying the following recurrence relation
\begin{equation*}
x_{1,0} = \dfrac{1-r}{\alpha_1 + \beta +1}, \quad x_{0,1} = \dfrac{r}{\alpha_2 + \beta +1},
\end{equation*}
\begin{equation} \label{rekurziox(w_1,w_2)}
x_{w_1,w_2} = \dfrac{\left(\alpha_1 \left( w_1-1 \right) + \beta_1\right)x_{w_1-1,w_2} + \left(\alpha_2 \left(w_2-1 \right) + \beta_2\right)x_{w_1,w_2-1}}
{\alpha_1 w_1 + \alpha_2 w_2 + \beta +1}
\end{equation}
if $1<w_1+w_2$.
\end{lem}

\begin{pf}
From the recurrence formula (\ref{rekurzio_x(d,w)}), we have

\begin{equation*}
x_{w_1,w_2}=\sum_{d=1}^{\left(N-1\right)w_1+w_2}x_{d,w_1,w_2}=\sum_{d}x_{d,w_1,w_2}=
\end{equation*}

\begin{equation*}
=\dfrac{1}{\alpha_1 w_1 + \alpha_2 w_2 + \beta +1} \times \left[ \alpha_{11} \left( w_1-1 \right) \sum_{d} x_{d-1,w_1-1,w_2} + \alpha_{12} \left(  w_1-1\right)\sum_{d}x_{d,w_1-1,w_2} +\right.
\end{equation*}
\begin{equation*}
+\alpha_2\left(w_2-1\right)\sum_{d}x_{d,w_1,w_2-1}+\left.\beta_1\sum_{d} x_{d-\left( N-1 \right),w_1-1,w_2} + \beta_2\sum_{d}x_{d-1,w_1,w_2-1} \right]=
\end{equation*}

\begin{equation*}
=\dfrac{\left(\alpha_1 \left( w_1-1 \right) + \beta_1\right)x_{w_1-1,w_2} + \left(\alpha_2 \left(w_2-1 \right) + \beta_2\right)x_{w_1,w_2-1}}{\alpha_1 w_1 + \alpha_2 w_2 + \beta +1}.
\end{equation*}
Above we used that $x_{d,w_1,w_2}=0$ if $d\notin\{1,2,\dots,\left(N-1\right)w_1+w_2\}$.

\qed
\end{pf}

\begin{thm} \label{thweight}
Assume $0<p<1$, $0<q<1$, $0<r<1$. 
Let $w_1$ be fixed.
Then, as $w_2\rightarrow\infty$,
\begin{equation}\label{w1fixthm}
x_{w_1,w_2}\sim C(w_1)w_2^{-\left(1+\frac{\beta_1+1}{\alpha_2}\right)},
\end{equation}
where
\begin{equation}\label{Cform}
C(w_1)=\dfrac{r}{\alpha_2}\dfrac{1}{w_1!}\dfrac{\varGamma\left(w_1+\dfrac{\beta_1}{\alpha_1}\right)}{\varGamma\left(\dfrac{\beta_1}{\alpha_1}\right)}\dfrac{\varGamma\left(1+\dfrac{\beta+1}{\alpha_2}\right)}{\varGamma\left(1+\dfrac{\beta_2}{\alpha_2}\right)}.
\end{equation}
Here $\varGamma$ denotes the Gamma function.
Similar result is true when $w_2$ is fixed and $w_1\rightarrow\infty$.
\end{thm}

\begin{pf}
Throughout the proof we shall use the following two facts for the $\Gamma$-function.
\begin{equation}\label{prudformula}
\sum_{i=0}^n\dfrac{\varGamma(i+a)}{\varGamma(i+b)}=\dfrac{1}{a-b+1}\left[\dfrac{\varGamma(n+a+1)}{\varGamma(n+b)}-\dfrac{\varGamma(a)}{\varGamma(b-1)}\right],
\end{equation}
see \cite{Prud}.
Applying Stirling's formula, we have
\begin{equation}\label{GSTR}
\dfrac{\varGamma\left(n+a\right)}{\varGamma\left(n+b\right)}\sim n^{-(b-a)}.
\end{equation}

We shall use mathematical induction.
As we have to handle the two-dimensional array $x_{w_1,w_2}$ satisfying equation \eqref{rekurziox(w_1,w_2)}, our formulae will be quite long. 
Therefore we divide the induction procedure into several steps.

{\bf (i) First consider the cases when one of the weights is equal to $0$.}
 Let $w_1=0$, then by Lemma \ref{xdw}, we obtain

\begin{equation}\label{x01}
x_{0,1}=\dfrac{r}{\alpha_2+\beta+1}>0,
\end{equation}

\begin{equation}\label{rek0l}
x_{0,l}=\dfrac{1}{l\alpha_2+\beta+1}\left((l-1)\alpha_2+\beta_2\right)x_{0,l-1},\;\;\;\;\;\;\; l>1.
\end{equation}

Using (\ref{rek0l}), we have

\begin{equation*}
x_{0,l}=x_{0,1}\prod_{j=2}^{l}\dfrac{\alpha_2(j-1)+\beta_2}{\alpha_2j+\beta+1}=
\dfrac{r}{\alpha_2+\beta+1} \cdot \dfrac{\alpha_2+\beta_2}{2\alpha_2+\beta+1} \cdot \dfrac{2\alpha_2+\beta_2}{3\alpha_2+\beta+1}\dots\dfrac{(l-1)\alpha_2+\beta_2}{l\alpha_2+\beta+1}=
\end{equation*}

\begin{equation*}
=\dfrac{r}{l\alpha_2+\beta+1}\prod_{j=1}^{l-1}\dfrac{\frac{\beta_2}{\alpha_2}+j}{\frac{\beta+1}{\alpha_2}+j}=\dfrac{r}{l\alpha_2+\beta+1}\dfrac{\varGamma\left(1+\frac{\beta+1}{\alpha_2}\right)}{\varGamma\left(1+\frac{\beta_2}{\alpha_2}\right)}\dfrac{\varGamma\left(l+\frac{\beta_2}{\alpha_2}\right)}{\varGamma\left(l+\frac{\beta+1}{\alpha_2}\right)}=
\end{equation*}

\begin{equation}\label{gammax0l}
=\dfrac{r}{\alpha_2}\dfrac{\varGamma\left(1+\frac{\beta+1}{\alpha_2}\right)}{\varGamma\left(1+\frac{\beta_2}{\alpha_2}\right)}\dfrac{\varGamma\left(l+\frac{\beta_2}{\alpha_2}\right)}{\varGamma\left(l+\frac{\alpha_2+\beta+1}{\alpha_2}\right)} = C(0) \dfrac{\varGamma\left(l+\frac{\beta_2}{\alpha_2}\right)}{\varGamma\left(l+\frac{\alpha_2+\beta+1}{\alpha_2}\right)}.
\end{equation}

Applying (\ref{GSTR}), we obtain our statement for $w_1=0$ that
\begin{equation}\label{x0lasmp}
x_{0,l}= C(0) \dfrac{\varGamma\left(l+\frac{\beta_2}{\alpha_2}\right)}{\varGamma\left(l+\frac{\alpha_2+\beta+1}{\alpha_2}\right)}\sim C(0) l^{-\left(1+\frac{\beta_1+1}{\alpha_2}\right)}
\end{equation}
as $l\rightarrow\infty$.

When $w_2=0$, we can obtain similarly the following

\begin{equation}\label{x10}
x_{1,0}=\dfrac{1-r}{\alpha_1+\beta+1}>0,
\end{equation}

\begin{equation}\label{rekk0}
x_{k,0}=\dfrac{1}{k\alpha_1+\beta+1}\left((k-1)\alpha_1+\beta_1\right)x_{k-1,0},\;\;\;\;\;\;\; k>1,
\end{equation}
and
\begin{equation}\label{xk0asmp}
x_{k,0}\sim C_0 k^{-\left(1+\frac{\beta_2+1}{\alpha_1}\right)}
\end{equation}
as $k\rightarrow\infty$, where $C_0=\dfrac{(1-r)\varGamma\left(1+\frac{\beta+1}{\alpha_1}\right)}{\alpha_1\varGamma\left(1+\frac{\beta_1}{\alpha_1}\right)}$.

{\bf (ii) Let us consider the case, when $w_1=1$.}
Analysing this particular case is very useful in order to find the general formulae which will be given in the next stage.
We have already seen that statement (\ref{w1fixthm}) is true if $w_1=0$ or $w_2=0$ (cf. (\ref{x0lasmp}) and \eqref{x10}).
 Applying several times (\ref{rekurziox(w_1,w_2)}), we obtain

\begin{equation*}
x_{1,l}=\dfrac{\beta_1x_{0,l} + \left(\alpha_2 \left(l-1 \right) + \beta_2\right)x_{1,l-1}}{\alpha_1 + \alpha_2 l + \beta +1}=
\end{equation*}

\begin{equation*}
=\sum_{i=1}^{l}\dfrac{\beta_1\prod_{j=i}^{l-1}\left(j\alpha_2+\beta_2\right)}{\prod_{j=i}^{l}\left(\alpha_1+j\alpha_2+\beta+1\right)}x_{0,i}+\dfrac{\beta_2\prod_{j=1}^{l-1}\left(j\alpha_2+\beta_2\right)}{\prod_{j=1}^{l}\left(\alpha_1+j\alpha_2+\beta+1\right)}x_{1,0}=
\end{equation*}

\begin{equation}\label{proof1}
=\sum_{i=1}^{l}b_{0,i}^{(l)}x_{0,i}+b_{1,0}^{(l)}x_{1,0},
\end{equation}
where

\begin{equation}\label{b0i}
b_{0,i}^{(l)}=\dfrac{\beta_1}{\alpha_2}\dfrac{\varGamma\left(l+\frac{\beta_2}{\alpha_2}\right)}{\varGamma\left(l+1+\frac{\alpha_1+\beta+1}{\alpha_2}\right)}\dfrac{\varGamma\left(i+\frac{\alpha_1+\beta+1}{\alpha_2}\right)}{\varGamma\left(i+\frac{\beta_2}{\alpha_2}\right)},\;\;\;\;\; 1\leq i\leq l,
\end{equation}
and

\begin{equation}\label{b10}
b_{1,0}^{(l)}=\dfrac{\beta_2}{\alpha_2}\dfrac{\varGamma\left(1+\frac{\alpha_1+\beta+1}{\alpha_2}\right)}{\varGamma\left(1+\frac{\beta_2}{\alpha_2}\right)}\dfrac{\varGamma\left(l+\frac{\beta_2}{\alpha_2}\right)}{\varGamma\left(l+1+\frac{\alpha_1+\beta+1}{\alpha_2}\right)}.
\end{equation}

Let
\begin{equation}\label{ACform0}
D(1)=\dfrac{C(0)\beta_1}{\alpha_2}\dfrac{\varGamma\left(l+\dfrac{\beta_2}{\alpha_2}\right)}{\varGamma\left(l+1+\dfrac{\alpha_1+\beta+1}{\alpha_2}\right)},
\end{equation}
where $C(0)$ is from (\ref{Cform}).

Now using (\ref{gammax0l}), (\ref{prudformula}), (\ref{GSTR}) for the first term in (\ref{proof1}), we obtain
\begin{equation*}
\sum_{i=1}^{l}b_{0,i}^{(l)}x_{0,i}= D(1) \sum_{i=1}^l\dfrac{\varGamma\left(i+\dfrac{\alpha_1+\beta+1}{\alpha_2}\right)}{\varGamma\left(i+\dfrac{\alpha_2+\beta+1}{\alpha_2}\right)}= D(1) \sum_{t=0}^{l-1}\dfrac{\varGamma\left(t+1+\dfrac{\alpha_1+\beta+1}{\alpha_2}\right)}{\varGamma\left(t+1+\dfrac{\alpha_2+\beta+1}{\alpha_2}\right)}=
\end{equation*}
\begin{equation*}
=\dfrac{D(1)}{\dfrac{\alpha_1+\beta+1}{\alpha_2}-\dfrac{\alpha_2+\beta+1}{\alpha_2}+1}\left[\dfrac{\varGamma\left(l+1+\dfrac{\alpha_1+\beta+1}{\alpha_2}\right)}{\varGamma\left(l+\dfrac{\alpha_2+\beta+1}{\alpha_2}\right)}-\dfrac{\varGamma\left(1+\dfrac{\alpha_1+\beta+1}{\alpha_2}\right)}{\varGamma\left(\dfrac{\alpha_2+\beta+1}{\alpha_2}\right)}\right]=
\end{equation*}
\begin{multline*}
=\dfrac{r}{\alpha_2}\dfrac{\beta_1}{\alpha_1}\dfrac{\varGamma\left(1+\dfrac{\beta+1}{\alpha_2}\right)}{\varGamma\left(1+\dfrac{\beta_2}{\alpha_2}\right)}
\left[\dfrac{\varGamma\left(l+\dfrac{\beta_2}{\alpha_2}\right)}{\varGamma\left(l+\dfrac{\alpha_2+\beta+1}{\alpha_2}\right)} +\right. \\ \left.-\dfrac{\varGamma\left(l+\dfrac{\beta_2}{\alpha_2}\right)}{\varGamma\left(l+1+\dfrac{\alpha_1+\beta+1}{\alpha_2}\right)}\dfrac{\varGamma\left(1+\dfrac{\alpha_1+\beta+1}{\alpha_2}\right)}{\varGamma\left(\dfrac{\alpha_2+\beta+1}{\alpha_2}\right)}\right]\sim
\end{multline*}
\begin{equation*}
\sim\dfrac{r}{\alpha_2}\dfrac{\beta_1}{\alpha_1}\dfrac{\varGamma\left(1+\dfrac{\beta+1}{\alpha_2}\right)}{\varGamma\left(1+\dfrac{\beta_2}{\alpha_2}\right)}
\dfrac{\varGamma\left(l+\dfrac{\beta_2}{\alpha_2}\right)}{\varGamma\left(l+\dfrac{\alpha_2+\beta+1}{\alpha_2}\right)}=C(1)\dfrac{\varGamma\left(l+\dfrac{\beta_2}{\alpha_2}\right)}{\varGamma\left(l+\dfrac{\alpha_2+\beta+1}{\alpha_2}\right)}
\end{equation*}
as $l\rightarrow\infty$.
Above we neglected the lower order term. Now using (\ref{GSTR}), we have
%
%
\begin{equation}\label{b0ix1iasmp}
\sum_{i=1}^{l}b_{0,i}^{(l)}x_{0,i}\sim C(1) l^{-\left(1+\frac{\beta_1+1}{\alpha_2}\right)}
\end{equation}
as $l\rightarrow\infty$.

Now, turn to the second term in \eqref{proof1}.
Applying (\ref{GSTR}) for (\ref{b10}), we obtain
\begin{equation}\label{b10asmp}
b_{1,0}^{(l)}\sim\dfrac{\beta_2}{\alpha_2}\dfrac{\varGamma\left(1+\frac{\alpha_1+\beta+1}{\alpha_2}\right)}{\varGamma\left(1+\frac{\beta_2}{\alpha_2}\right)}l^{-\left(1+\frac{\alpha_1+\beta_1+1}{\alpha_2}\right)}
\end{equation}
as $l\rightarrow\infty$.

Now substitute (\ref{x10}), (\ref{b0ix1iasmp}), (\ref{b10asmp}) into (\ref{proof1}), we obtain the statement for $w_1=1$
\begin{equation*}
x_{1,l}\sim C(1) l^{-\left(1+\frac{\beta_1+1}{\alpha_2}\right)} + \dfrac{\beta_2}{\alpha_2}\dfrac{\varGamma\left(1+\frac{\alpha_1+\beta+1}{\alpha_2}\right)}{\varGamma\left(1+\frac{\beta_2}{\alpha_2}\right)}l^{-\left(1+\frac{\alpha_1+\beta+1}{\alpha_2}\right)}\dfrac{1-r}{\alpha_1+\beta+1}\sim
\end{equation*}
\begin{equation*}
\sim C(1) l^{-\left(1+\frac{\beta_1+1}{\alpha_2}\right)}
\end{equation*}
as $l\rightarrow\infty$.

{\bf (iii) Now consider the general induction step.}
Suppose that the statement is true for all central weights being less than $w_1$. 
That is, the following is true
\begin{equation}\label{statement}
x_{k,l}\sim C(k) \dfrac{\varGamma\left(l+\frac{\beta_2}{\alpha_2}\right)}{\varGamma\left(l+\frac{\alpha_2+\beta+1}{\alpha_2}\right)}
\sim C(k)l^{-\left(1+\frac{\beta_1+1}{\alpha_2}\right)},\;\;\;\;\; 0\leq k\leq w_1-1,
\end{equation}
as $l\rightarrow\infty$.
Then, from (\ref{rekurziox(w_1,w_2)}), we get
\begin{equation*}
x_{w_1,l}=\dfrac{(w_1-1)\alpha_1+\beta_1}{\alpha_2} \sum_{i=1}^{l}\dfrac{\prod_{j=i}^{l-1}\left(j+\frac{\beta_2}{\alpha_2}\right)}{\prod_{j=i}^{l}\left(j+\frac{w_1\alpha_1+\beta+1}{\alpha_2}\right)}x_{w_1-1,i}+\dfrac{\beta_2}{\alpha_2}\dfrac{\prod_{j=1}^{l-1}\left(j+\frac{\beta_2}{\alpha_2}\right)}{\prod_{j=1}^{l}\left(j+\frac{w_1\alpha_1+\beta+1}{\alpha_2}\right)}x_{w_1,0}=
\end{equation*}
\begin{equation}\label{xw1l}
=\sum_{i=1}^lb_{w_1-1,i}^{(l)}x_{w_1-1,i}+b_{w_1,0}^{(l)}x_{w_1,0},
\end{equation}
where
\begin{equation}\label{bw11i}
b_{w_1-1,i}^{(l)}=\dfrac{(w_1-1)\alpha_1+\beta_1}{\alpha_2}\dfrac{\varGamma\left(l+\frac{\beta_2}{\alpha_2}\right)}{\varGamma\left(l+1+\frac{w_1\alpha_1+\beta+1}{\alpha_2}\right)}\dfrac{\varGamma\left(i+\frac{w_1\alpha_1+\beta+1}{\alpha_2}\right)}{\varGamma\left(i+\frac{\beta_2}{\alpha_2}\right)},\;\;\;\;\; 1\leq i\leq l,
\end{equation}
and

\begin{equation}\label{bw10}
b_{w_1,0}^{(l)}=\dfrac{\varGamma\left(1+\frac{w_1\alpha_1+\beta+1}{\alpha_2}\right)}{\varGamma\left(\frac{\beta_2}{\alpha_2}\right)}\dfrac{\varGamma\left(l+\frac{\beta_2}{\alpha_2}\right)}{\varGamma\left(l+1+\frac{w_1\alpha_1+\beta+1}{\alpha_2}\right)}.
\end{equation}

%
%
Let
\begin{equation}\label{ACform}
D(w_1)=\dfrac{C(w_1-1)\left((w_1-1)\alpha_1+\beta_1\right)}{\alpha_2}\dfrac{\varGamma\left(l+\dfrac{\beta_2}{\alpha_2}\right)}{\varGamma\left(l+1+\dfrac{w_1\alpha_1+\beta+1}{\alpha_2}\right)},
\end{equation}
where $C(w_1-1)$ is from (\ref{Cform}).

Using (\ref{statement}), (\ref{prudformula}), (\ref{GSTR}) and (\ref{ACform}) for the first term in (\ref{xw1l}), we obtain
\begin{multline*}
\sum_{i=1}^{l}b_{w_1-1,i}^{(l)}x_{w_1-1,i} \sim D(w_1) \sum_{i=1}^l\dfrac{\varGamma\left(i+\dfrac{w_1\alpha_1+\beta+1}{\alpha_2}\right)}{\varGamma\left(i+\dfrac{\alpha_2+\beta+1}{\alpha_2}\right)}=\\=D(w_1) \sum_{t=0}^{l-1}\dfrac{\varGamma\left(t+1+\dfrac{w_1\alpha_1+\beta+1}{\alpha_2}\right)}{\varGamma\left(t+1+\dfrac{\alpha_2+\beta+1}{\alpha_2}\right)}=
\end{multline*}
\begin{multline*}
=\dfrac{D(w_1)}{\dfrac{w_1\alpha_1+\beta+1}{\alpha_2}-\dfrac{\alpha_2+\beta+1}{\alpha_2}+1}\left[\dfrac{\varGamma\left(l+1+\dfrac{w_1\alpha_1+\beta+1}{\alpha_2}\right)}{\varGamma\left(l+\dfrac{\alpha_2+\beta+1}{\alpha_2}\right)}+\right.\\\left.-\dfrac{\varGamma\left(1+\dfrac{w_1\alpha_1+\beta+1}{\alpha_2}\right)}{\varGamma\left(\dfrac{\alpha_2+\beta+1}{\alpha_2}\right)}\right]=
\end{multline*}
\begin{equation*}
=\dfrac{r}{\alpha_2}\dfrac{1}{w_1!}\dfrac{\varGamma\left(w_1+\dfrac{\beta_1}{\alpha_1}\right)}{\varGamma\left(\dfrac{\beta_1}{\alpha_1}\right)}\dfrac{\varGamma\left(1+\dfrac{\beta+1}{\alpha_2}\right)}{\varGamma\left(1+\dfrac{\beta_2}{\alpha_2}\right)}
\left[\dfrac{\varGamma\left(l+\frac{\beta_2}{\alpha_2}\right)}{\varGamma\left(l+\dfrac{\alpha_2+\beta+1}{\alpha_2}\right)}+\right.
\end{equation*}
\begin{equation*}
 -\left. \dfrac{\varGamma\left(l+\dfrac{\beta_2}{\alpha_2}\right)}{\varGamma\left(l+1+\dfrac{w_1\alpha_1+\beta+1}{\alpha_2}\right)}\dfrac{\varGamma\left(1+\dfrac{w_1\alpha_1+\beta+1}{\alpha_2}\right)}{\varGamma\left(\dfrac{\alpha_2+\beta+1}{\alpha_2}\right)}\right]\sim C(w_1)\dfrac{\varGamma\left(l+\frac{\beta_2}{\alpha_2}\right)}{\varGamma\left(l+\dfrac{\alpha_2+\beta+1}{\alpha_2}\right)}
\end{equation*}
as $l\rightarrow\infty$.
Above we neglected the lower order term. Using (\ref{GSTR}), we have
%
%
\begin{equation}\label{bw11ixw11iasmp}
\sum_{i=1}^{l}b_{w_1-1,i}^{(l)}x_{w_1-1,i}\sim C(w_1) l^{-\left(1+\frac{\beta_1+1}{\alpha_2}\right)}
\end{equation}
as $l\rightarrow\infty$.

Applying (\ref{GSTR}) for (\ref{bw10}), we obtain
\begin{equation}\label{bw10asmp}
b_{w_1,0}^{(l)}\sim\dfrac{\beta_2}{\alpha_2}\dfrac{\varGamma\left(1+\frac{w_1\alpha_1+\beta+1}{\alpha_2}\right)}{\varGamma\left(1+\frac{\beta_2}{\alpha_2}\right)}l^{-\left(1+\frac{w_1\alpha_1+\beta_1+1}{\alpha_2}\right)}
\end{equation}
as $l\rightarrow\infty$.

Now substitute (\ref{bw11ixw11iasmp}), (\ref{bw10asmp}) and (\ref{xk0asmp})  into (\ref{xw1l}), and we obtain our statement
\begin{equation*}
x_{w_1,l}\sim C(w_1) l^{-\left(1+\frac{\beta_1+1}{\alpha_2}\right)} +\dfrac{\beta_2}{\alpha_2}\dfrac{\varGamma\left(1+\frac{w_1\alpha_1+\beta+1}{\alpha_2}\right)}{\varGamma\left(1+\frac{\beta_2}{\alpha_2}\right)}l^{-\left(1+\frac{w_1\alpha_1+\beta_1+1}{\alpha_2}\right)} C_0 w_1^{-\left(1+\frac{\beta_2+1}{\alpha_1}\right)} \sim
\end{equation*}
\begin{equation*}
\sim C(w_1) l^{-\left(1+\frac{\beta_1+1}{\alpha_2}\right)}
\end{equation*}
as $l\rightarrow\infty$.
\qed
\end{pf}
%
\section{Appendix}
\setcounter{equation}{0}
We use the following results on discrete time submartingales.
Let $\{Z_n, \FD_n \} $ be a submartingale.
Its Doob-Meyer decomposition is $Z_n = M_n + A_n$, where $\{M_n, \FD_n \}$ is a martingale and $\{A_n, \FD_n \}$ is an increasing predictable process.
Here, up to an additive constant, 
$$
A_n= \E Z_1 + \sum_{i=2}^n (\E(Z_i | \FD_{i-1}) - Z_{i-1} ).
$$
Now, $\{M^2_n, \FD_n \}$ is again a submartingale.
Let 
$$
M^2_n = Y_n + B_n
$$
be the Doob-Meyer decomposition of $M^2_n$.
Here, up to an additive constant, 
$$
B_n= \sum_{i=2}^n \D^2(Z_i | \FD_{i-1}) = \sum_{i=2}^n \E\left\{ (Z_i-\E(Z_i | \FD_{i-1}))^2 | \FD_{i-1} \right\}.
$$
\begin{lem} \label{neveu}
Let $M_1=0$.
On the set $\{ B_\infty < \infty \}$ the martingale $M_n$ almost surely converges to a finite limit.
Moreover, $M_n = \o (B_n^{1/2} \log{B_n} )$ almost surely on the set $\{ B_n \to \infty \}$.
\newline
Let $\{Z_n, \FD_n \} $ be a square integrable non-negative submartingale.
If  $B_n^{1/2} \log{B_n} = \O (A_n)$, then $Z_n \sim A_n$ as $n\to\infty$, almost surely on the set $\{ A_n \to \infty \}$.
\end{lem}
The first part of the lemma is contained in Propositions VII-2-3 and VII-2-4 of \cite{Neveu}.
The second part is a consequence of the first part and it can be found in Proposition 2.3 of \cite{backhausz}.

\bigskip
\noindent
{\bf Acknowledgements.} The  authors are grateful to the referee and the editor for the careful reading of the paper and for the valuable suggestions.


\begin{thebibliography}{99}
\bibitem{backhausz}
Backhausz, \'A., \textit{Analysis of random graphs with methods of martingale theory.} PhD thesis, E\"otv\"os Lor\'and University, Budapest, 2012.
%
\bibitem{BaMo1}
Backhausz, {\'A.,} M\'ori, T.\,F.,  \textit{A random graph model based on 3-interactions}. Ann. Univ. Sci. Budapest. Sect. Comput., {\bf 36} (2012), 41--52.
%
\bibitem{BaMo2}
Backhausz, {\'A.,} M\'ori, T.\,F., \textit{Weights and degrees in a random graph model based on 3-interactions}. Acta Math. Hungar., {\bf 143} no.1 (2014), 23--43.
%
\bibitem{barabasiBook}
Barab\'asi, A.\,L., \textit{Network Science}. Cambridge University Press, 2016.
%
\bibitem{barabasi}
Barab\'asi, A.\,L., Albert, R., \textit{Emergence of scaling in random networks}. Science, {\bf 286} (1999), 509--512.
%
\bibitem{bollobas}
Bollob\'as, B., Riordan, O.\,M., Spencer, J., Tusn\'ady, G., \textit{The degree sequence of a scale-free random graph process}. Random Structures Algorithms, {\bf 18} (2001), 279--290.
%
\bibitem{bollobasRiordan}
Bollob\'as, B., Riordan, O.\,M., \textit{Mathematical results on scale-free random graphs}. In: Handbook of Graphs and Networks; From the Genome to the Internet, pp. 1--34, Wiley, 2003.
%
\bibitem{comp}
Bosworth, S., Kabay, M.\,E., Whyne, E. (editors), \textit{Computer Security Handbook}.
John Wiley \& Sons, 2014.
%
\bibitem{broido}
Broido, A.\,D.,  Clauset, A., \textit{Scale-free networks are rare}. arXiv:1801.03400, 2018.
%
\bibitem{buckleyOsthus}
Buckley, P.\,G., Osthus, D., \textit{ Popularity based random graph models leading to a scale-free degree sequence}. Discrete Mathematics, {\bf 282} no.1-6 (2004), 53--63.
%
\bibitem{Chung-Lu}
Chung, F., Lu, L., \textit{Complex Graphs and Networks}. AMS and CBMS, 2006.
%
\bibitem{cooper}
Cooper, C., Frieze, A., \textit{A general model of web graphs}. Random Structures Algorithms, {\bf 22} (2003), 311--335.
%
\bibitem{durrett}
Durrett, R., \textit{Random graph dynamics}. Cambridge University Press, Cambridge, 2007.
%
\bibitem{erRe1}
Erd\H{o}s, P., R\'enyi, A.,  \textit{On random graphs, I}. Publicationes Mathematicae (Debrecen), {\bf 6} (1959), 290--297.
%
\bibitem{erRe2}
Erd\H{o}s, P., R\'enyi, A., \textit{On the evolution of random graphs}. Publ. Math. Inst. Hung. Acad. Sci., {\bf 5} (1960), 17--61. 
%
%
\bibitem{FIPB2}
Fazekas, I., Porv\'azsnyik, B., \textit{Scale-free property for degrees and weights in an $N$-interactions random graph model}.
Journal of Mathematical Sciences, {\bf 214} no.1 (2016), 69--82.
%
\bibitem{FIPB3}
Fazekas, I., Porv\'azsnyik, B., \textit{Limit theorems for the weights and the degrees in an $N$-interactions random graph model}.
Open Mathematics, {\bf 14} no.1 (2016), 414--424.
%
\bibitem{Faz} Fazekas, I., Nosz\'aly, Cs., Perecs\'enyi, A., 
\textit{A population evolution model and its applications to random networks}. arXiv:1604.01579, 2016.
%
\bibitem{gil}
Gilbert, E.\,N., \textit{Random graphs}. The Annals of Mathematical Statistics, {\bf 30} (1959), 1141--1144. 
%
\bibitem{holmeKim}
Holme, P., Kim, B.\,J., \textit{Growing scale-free networks with tunable clustering}. Phys. Rev. E {\bf 65} no.2, 026107 (2002).
%
\bibitem{Mori}
M\'ori, T.\,F., \textit{The maximum degree of the Barab\'asi-Albert random tree}. Combinatorics, Probability and Computing, {\bf 14} no.3 (2005), 339--348.
%
\bibitem{Neveu}
Neveu, J., \textit{Discrete-parameter martingales.} North-Holland, Amsterdam, 1975.
%
\bibitem{osroumova} Ostroumova, L.,  Ryabchenko, A., Samosvat, E.,
\textit{Generalized preferential attachment: tunable power-law degree distribution and clustering coefficient}.
Algorithms and Models for the Web Graph, Lecture Notes in Computer Science, volume 8305, pp. 185--202, Springer, 2013.
%
\bibitem{Prud}
Prudnikov, A.\,P.; Brychkov, Yu.\,A.; Marichev, O.\,I., \textit{Integrals and series.} Gordon \& Breach Science Publishers, New York, 1986.
%
%
\bibitem{Simon} 
Simon, H.\,A.,  \textit{On a class of skew distribution functions}. Biometrika, {\bf 42} no.3-4 (1955), 425--440.
%
%
\bibitem{hofstad}
van der Hofstad, R., \textit{Random Graphs and Complex Networks.} Cambridge University Press, 2017.
%
\bibitem{Zhou}
Zhou, T., Yan, G., Wang, B.-H., \textit{Maximal planar networks with large clustering coefficient and power-law degree distribution}. Phys. Rev. E {\bf 71} no.4 (2005), 046141.
Erratum Phys. Rev. E 72, 029905 (2005).
%
\bibitem{Yule}
Yule, G.\,U., \textit{A mathematical theory of evolution, based on the conclusions of Dr. J. C. Willis, F. R. S.} Phil. Trans. Roy. Soc. London, B, {\bf 213} (1925), 21--87.
%
\end{thebibliography}
\end{document}